\documentclass[11pt]{amsart}

\usepackage{epigamath}

\usepackage[english]{babel}

\usepackage[top=27 mm, bottom=20 mm, left=25 mm, right=25mm]{geometry}
\usepackage[all]{xy}
\usepackage{multirow}
\usepackage{enumitem}
\usepackage{mathtools,stackengine}
\usepackage{pdflscape}
\usepackage{arydshln}

\theoremstyle{definition}
\newtheorem{Definition}{Definition}[section]

\newtheorem{Example}[Definition]{Example}

\newtheorem{Remark}[Definition]{Remark}
\newtheorem{Remark/Notation}[Definition]{Remark/Notation}

\newtheorem*{Algorithmx}{Algorithm}

\newtheorem{Question}[Definition]{Question}

\newtheorem*{Notationx}{Notation}

\theoremstyle{plain}
\newtheorem{Theorem}[Definition]{Theorem}

\newtheorem*{Main Theoremx}{Main Theorem}
\newtheorem{Proposition}[Definition]{Proposition}
\newtheorem{Observation}[Definition]{Observation}

\newtheorem{Lemma}[Definition]{Lemma}

\newtheoremstyle{voiditstyle}{3pt}{3pt}{\itshape}{\parindent}%
{\bfseries}{.}{ }{\thmnote{#3}}%
\theoremstyle{voiditstyle}

\newtheoremstyle{voidromstyle}{3pt}{3pt}{\rm}{\parindent}%
{\bfseries}{.}{ }{\thmnote{#3}}%
\theoremstyle{voidromstyle}

\newcommand{\prf}{\par\noindent{\sc Proof.}\quad}

\newcommand{\cal}{\mathcal}

\newcommand{\bbZ}{{\mathbb{Z}}}
\newcommand{\bbP}{{\mathbb{P}}}

\newcommand{\Spec}{\operatorname{Spec}}

\newcommand{\Pic}{\operatorname{Pic}}

\newcommand{\Char}{\operatorname{char}}

\newcommand{\bsm}{\left(\begin{smallmatrix}}
\newcommand{\esm}{\end{smallmatrix}\right)}

\newcommand{\beq}{\begin{equation}}
\newcommand{\eeq}{\end{equation}}

\numberwithin{equation}{section}

\EpigaVolumeYear{5}{2021} \EpigaArticleNr{17}
\ReceivedOn{January 4, 2021}
\InFinalFormOn{July 6, 2021}
\AcceptedOn{October 25, 2021}

\title{Which rational double points occur on del Pezzo surfaces?}
\titlemark{Which rational double points occur on del Pezzo surfaces?}

\author{Claudia Stadlmayr}
\address{TU M\"unchen, Zentrum Mathematik - M11, Boltzmannstra{\ss}e 3, 85748 Garching bei M\"unchen, Germany}
\curraddr{}
\email{claudia.stadlmayr@ma.tum.de}

\authormark{C.~Stadlmayr}

\AbstractInEnglish{We determine all configurations of rational double points that occur on RDP del Pezzo surfaces of arbitrary degree and Picard rank over an algebraically closed field $k$ of arbitrary characteristic $\Char(k)=p \geq 0$, generalizing classical work of Du Val to positive characteristic. Moreover, we give simplified equations for all RDP del Pezzo surfaces of degree $1$ containing non-taut rational double points.}



\TitleInFrench{Quels points doubles rationels trouve-t-on sur les surfaces de del Pezzo ?}

\AbstractInFrench{On d\'etermine toutes les configurations de points doubles rationnels apparaissant sur les surfaces
de del Pezzo n'ayant que des points doubles rationnels, de degr\'e et rang du groupe de Picard arbitraires, sur des
corps alg\'ebriquement clos $k$ de caract\'eristique quelconque $\Char(k)=p \geq 0$. Ceci g\'en\'eralise le travail classique de
Du Val \`a la caract\'eristique positive. De plus, on donne des \'equations simplifi\'ees pour les surfaces de del Pezzo de degr\'e $1$ n'ayant que des points doubles rationnels dont au moins une singularit\'e non \emph{taut}.}

\acknowledgement{Research of the author is funded by the DFG Sachbeihilfe LI 1906/5-1 ``Geometrie von rationalen Doppelpunkten''. The author gratefully acknowledges support by the doctoral program TopMath and the TUM Graduate School.}

\begin{document}

%%%%%%%%%%%%%%%%%%%%%%%%%%%%%%%
% Title page
%%%%%%%%%%%%%%%%%%%%%%%%%%%%%%%

%\removeabove{}
%\removebetween{}
%\removebelow{}

\maketitle

\begin{prelims}

\DisplayAbstractInEnglish

\bigskip

\DisplayKeyWords

\medskip

\DisplayMSCclass

\bigskip

\languagesection{Fran\c{c}ais}

\bigskip

\DisplayTitleInFrench

\medskip

\DisplayAbstractInFrench

\end{prelims}

%%%%%%%%%%%%%%%%%%%%%
% Table of Contents
%%%%%%%%%%%%%%%%%%%%%

\newpage

\setcounter{tocdepth}{1}

\tableofcontents

%%%%%%%%%%%%%%%%%%%%%
% Content begins here
%%%%%%%%%%%%%%%%%%%%%
\vspace*{-1cm}
\section{Introduction}

Throughout this article, we work over an algebraically closed field $k$ of characteristc $\Char(k) = p \geq 0$. Quite classically, del Pezzo surfaces admitting at worst rational double point singularities (also called \emph{RDP del Pezzo surfaces}, \emph{Gorenstein log del Pezzo surfaces} or \emph{Du Val del Pezzo surfaces}) first appeared as non-degenerate surfaces of degree $d$ in $\bbP^d$ which are not cones or projections of surfaces of minimal degree. A first natural question to ask is the following: 

\begin{Question} \label{Question central!}
Which rational double points occur on RDP del Pezzo surfaces?
\end{Question}

For $1 \leq d \leq 9$, we define the lattice $E_{9-d}$ as the orthogonal complement of the vector $(-3,1,\hdots,1)$ in the unimodular lattice ${\rm I}^{1,9-d}$ of signature $(1,9-d)$ defined by the matrix ${\rm diag}(1,-1,\hdots,-1)$. Over the complex numbers, the answer to Question \ref{Question central!} is old: 
\begin{itemize}
    \item If $d \geq 3$, a configuration $\Gamma$ of rational double points occurs on an RDP del Pezzo surface $X_d$ of degree $d$ if and only if the lattice $\Gamma'$ spanned by the irreducible components of the exceptional divisor of the minimal resolution $\widetilde{X}$ of $X$ embeds into the lattice $E_{9-d}$. If $d \geq 5$, this is elementary to check, while if $d = 4$ it was proven by Timms \cite{Timms} in 1928, and if $d = 3$ it was proven by Schl{\"a}fli \cite{Schlafli} in 1863.
    \item If $d = 2$, then $X_2$ is a double cover of $\bbP^2$ branched over a quartic curve. Simple singularities of plane quartics, and thus RDPs on $X_2$, were classified by Du Val \cite{DuValI-II-III} in 1934. It turns out that $\Gamma$ occurs on some $X_2$ if and only if $\Gamma'$ embeds into $E_7$ and $\Gamma$ is not of type $7A_1$.
    \item If $d = 1$, then $X_1$ is a double cover of a quadric cone in $\bbP^3$ branched over a sextic curve and the possible singularities of the branch locus have also been classified by Du Val \cite{DuValI-II-III}. It turns out that $\Gamma$ occurs on some $X_1$ if and only if $\Gamma'$ embeds into $E_8$ and $\Gamma \not \in \{D_4 + 4A_1,8A_1,7A_1\}$.
\end{itemize}
Besides these very classical sources, we refer the reader to \cite{Urabe} for a more modern treatment of the cases where $d \in \{1,2\}$. From the above discussion, we see that if we do not care about the degree of $X$, but only about whether $\Gamma$ occurs on \emph{some} RDP del Pezzo surface, then the lattice $E_8$ plays a central role and we note that the classification of root sublattices of $E_8$ is also quite classical and goes back to Dynkin \cite{DynkinSemisimple}. 

Thus, we are mainly interested in Question \ref{Question central!} in the case $p > 0$, even though our methods also recover Du Val's results over the complex numbers. While an answer to Question \ref{Question central!} may be known to the experts if $p \neq 2,3,5$, it becomes particularly subtle in small characteristics, where a rational double point is not necessarily \emph{taut}, that is, it is not necessarily uniquely determined by its dual resolution graph. Such non-taut rational double points can only occur in characteristic $2,3$, and $5$, and have first been studied and classified by Artin \cite{ArtinCoveringsRDP}: there are only finitely many rational double points with the same resolution graph and they are distinguished by a coindex that we call \emph{Artin coindex} (see \mbox{Table \ref{TableNonTaut-char235}} for a summary of all non-taut rational double points).

Very recently, there has been substantial progress on RDP del Pezzo surfaces in positive characteristic: 
\begin{itemize}
    \item In characteristic at least $5$, RDP del Pezzo surfaces of Picard rank $1$ have been classified by Lacini in \cite{Lacini}, generalizing work of Ye \cite{Ye} and Furushima \cite{Furushima} (see also \cite{MiyanishiZhang1} and \cite{MiyanishiZhang2}).
    \item In characteristic $2$ and $3$, Kawakami and Nagaoka \cite{KawakamiNagaokaNeu-PicardRank1} classify RDP del Pezzo surfaces of Picard rank $1$ and determine some, but not all, of the Artin coindices of the rational double points that occur. In \cite{KawakamiNagaokaNeu-Pathologies}, they also investigate in detail some interesting pathological examples in characteristic $2$ that will also appear as exceptional cases in the present article (see Proposition \ref{Prop RDP iff RDP deg 1}).
\end{itemize}
\noindent
In this article, instead of studying RDP del Pezzo surfaces of small Picard rank, which play a prominent role in the minimal model program (see \textit{e.g.} \cite[Lemma 2]{MiyanishiZhang2}), we want to approach the problem from a more classical angle and try to find a satisfying positive characteristic analogue of Du Val's work relating Question \ref{Question central!} to the lattice $E_8$.

Before stating our main result, let us fix some terminology. 
We say that a sum $\Gamma' = \sum_{i} \Gamma_{i,n_i}$ of root lattices (that is, $\Gamma_i \in \{A,D,E\}$ and $n_i$ is the number of simple roots) occurs on a weak del Pezzo surface $\widetilde{X}$, if it is isomorphic to the lattice spanned by all $(-2)$-curves on $\widetilde{X}$.

Then, Theorem \ref{Thm main} gives a complete answer to Question \ref{Question central!}.

\vspace{-0.3mm}
\begin{Theorem} \label{Thm main}
Let $\Gamma = \sum_i \Gamma_{i,n_i}^{k_i}$ be an RDP configuration with Artin coindices $k_i$ and let $\Gamma' = \sum_{i} \Gamma_{i,n_i}$ be the lattice spanned by the irreducible components of the exceptional divisor of its minimal resolution.
\begin{enumerate}
    \item \label{label Thm main 1} If $p \neq 2$, then the following are equivalent:
    \begin{itemize}
        \item $\Gamma$ occurs on an RDP del Pezzo surface.
		\item $\Gamma'$ occurs on a weak del Pezzo surface.
        \item $\Gamma'$ embeds into $E_8$ and $\Gamma' \not \in \{D_4 + 4A_1,8A_1,7A_1\}$.
    \end{itemize}
    \item \label{label Thm main 2} If $p = 2$, then the following are equivalent:
    \begin{itemize}
        \item $\Gamma'$ occurs on a weak del Pezzo surface.
        \item $\Gamma'$ embeds into $E_8$ and $\Gamma' \not \in \{2A_3 + 2A_1, A_3 + 4A_1, 6A_1\}$.
    \end{itemize}
    If all $\Gamma_{i,n_i}^{k_i}$ are taut, these statements are also equivalent to the following: 
    \begin{itemize}
    \item $\Gamma$ occurs on an RDP del Pezzo surface.
    \end{itemize}
    
    \item \label{label Thm main 3} If $p = 2$ and some $\Gamma_{i,n_i}^{k_i}$ is non-taut, then the following are equivalent:
    \begin{itemize}
        \item $\Gamma$ occurs on an RDP del Pezzo surface.
        \item $\Gamma$ occurs in Table \ref{TableD4D5-char2}, \ref{TableD6D7D8-char2}, or \ref{TableE6E7E8-char2}.
    \end{itemize}
    \item \label{label Thm main 4}
    Moreover, in Tables \ref{TableD4D5-char2}, \ref{TableD6D7D8-char2}, \ref{TableE6E7E8-char2}, \ref{TableE6E7E8-char3}, and \ref{TableE8-char5},
    we give equations for all RDP del Pezzo surfaces of degree $1$ containing a non-taut rational double point.
\end{enumerate}
\end{Theorem}

For the convenience of the reader, we list all possible $\Gamma'$ embedding into $E_8$ (see \cite[Table 11]{DynkinSemisimple}) and state whether the respective $\Gamma'$ occurs on a weak del Pezzo surface or not. For the possible Artin coindices of the corresponding RDP configurations in the non-taut case if $p=2$, we refer the reader to Tables \ref{TableD4D5-char2}, \ref{TableD6D7D8-char2}, and \ref{TableE6E7E8-char2} in the appendix (see also Remark \ref{Rem Artin coindices which do not occur}).

\vspace{2mm}

\begin{table}[h] \renewcommand{\arraystretch}{1.35}
 $$
 \begin{array}{|c||c|c|} 
 \hline
 
\multirow{2}{*}{$\Gamma' \hookrightarrow E_8$} 
 & \multicolumn{2}{c|}{\text{occurs if}} \\
 & {\text{$p \neq 2$}}
 & {\text{$p = 2$}}
 \\ \hline \hline

 A_1
 & \scriptstyle{\checkmark}
 & \scriptstyle{\checkmark}
 \\ \hline
 
  2A_1
 & \scriptstyle{\checkmark}
 & \scriptstyle{\checkmark}
 \\ \hline
 
  A_2
 & \scriptstyle{\checkmark}
 & \scriptstyle{\checkmark}
 \\ \hline
 
  3A_1
 & \scriptstyle{\checkmark}
 & \scriptstyle{\checkmark}
 \\ \hline
 
 A_2+A_1
 & \scriptstyle{\checkmark}
 & \scriptstyle{\checkmark}
 \\ \hline
 
  A_3
 & \scriptstyle{\checkmark}
 & \scriptstyle{\checkmark}
 \\ \hline
 
  4A_1
 & \scriptstyle{\checkmark}
 & \scriptstyle{\checkmark}
 \\ \hline
 
  A_2+2A_1
 & \scriptstyle{\checkmark}
 & \scriptstyle{\checkmark}
 \\ \hline
 
 2A_2
 & \scriptstyle{\checkmark}
 & \scriptstyle{\checkmark}
 \\ \hline
 
 A_3+A_1
 & \scriptstyle{\checkmark}
 & \scriptstyle{\checkmark}
 \\ \hline
 
  A_4
 & \scriptstyle{\checkmark}
 & \scriptstyle{\checkmark}
 \\ \hline
 
  D_4
 & \scriptstyle{\checkmark}
 & \scriptstyle{\checkmark}
 \\ \hline
 
  5A_1
 & \scriptstyle{\checkmark}
 & \scriptstyle{\checkmark}
 \\ \hline
 
  A_2+3A_1
 & \scriptstyle{\checkmark}
 & \scriptstyle{\checkmark}
 \\ \hline
 
 2A_2+A_1
 & \scriptstyle{\checkmark}
 & \scriptstyle{\checkmark}
 \\ \hline
 
  A_3+2A_1
 & \scriptstyle{\checkmark}
 & \scriptstyle{\checkmark}
 \\ \hline
 
  A_3+A_2
 & \scriptstyle{\checkmark}
 & \scriptstyle{\checkmark}
 \\ \hline
 
  A_4+A_1
 & \scriptstyle{\checkmark}
 & \scriptstyle{\checkmark}
 \\ \hline
 
   D_4+A_1
 & \scriptstyle{\checkmark}
 & \scriptstyle{\checkmark}
 \\ \hline
 
   A_5
 & \scriptstyle{\checkmark}
 & \scriptstyle{\checkmark}
 \\ \hline
 
  D_5
 & \scriptstyle{\checkmark}
 & \scriptstyle{\checkmark}
 \\ \hline

 6A_1
 & \scriptstyle{\checkmark}
 & \times
 \\ \hline
 
  A_2+4A_1
 & \scriptstyle{\checkmark}
 & \scriptstyle{\checkmark}
 \\ \hline

  2A_2+2A_1
 & \scriptstyle{\checkmark}
 & \scriptstyle{\checkmark}
 \\ \hline
 
  \end{array}
\hspace{1.5mm}
 \begin{array}{|c||c|c|} 
 \hline
 
\multirow{2}{*}{$\Gamma' \hookrightarrow E_8$} 
 & \multicolumn{2}{c|}{\text{occurs if}} \\
 & {\text{$p \neq 2$}}
 & {\text{$p = 2$}}
 \\ \hline \hline

  A_3+3A_1
 & \scriptstyle{\checkmark}
 & \scriptstyle{\checkmark}
 \\ \hline
 
   3A_2
 & \scriptstyle{\checkmark}
 & \scriptstyle{\checkmark}
 \\ \hline
 
  A_3+A_2+A_1
 & \scriptstyle{\checkmark}
 & \scriptstyle{\checkmark}
 \\ \hline

 A_4+2A_1
 & \scriptstyle{\checkmark}
 & \scriptstyle{\checkmark}
 \\ \hline
 
  D_4+2A_1
 & \scriptstyle{\checkmark}
 & \scriptstyle{\checkmark}
 \\ \hline

  2A_3
 & \scriptstyle{\checkmark}
 & \scriptstyle{\checkmark}
 \\ \hline
 
  A_4+A_2
 & \scriptstyle{\checkmark}
 & \scriptstyle{\checkmark}
 \\ \hline
 
 D_4+A_2
 & \scriptstyle{\checkmark}
 & \scriptstyle{\checkmark}
 \\ \hline
 
 A_5+A_1
 & \scriptstyle{\checkmark}
 & \scriptstyle{\checkmark}
 \\ \hline
 
  D_5+A_1
 & \scriptstyle{\checkmark}
 & \scriptstyle{\checkmark}
 \\ \hline
 
  A_6
 & \scriptstyle{\checkmark}
 & \scriptstyle{\checkmark}
 \\ \hline
 
  D_6
 & \scriptstyle{\checkmark}
 & \scriptstyle{\checkmark}
 \\ \hline
 
  E_6
 & \scriptstyle{\checkmark}
 & \scriptstyle{\checkmark}
 \\ \hline
 
 7A_1
 & \times
 & \scriptstyle{\checkmark}
 \\ \hline
 
  A_3+4A_1
 & \scriptstyle{\checkmark}
 & \times
 \\ \hline
 
  3A_2+A_1
 & \scriptstyle{\checkmark}
 & \scriptstyle{\checkmark}
 \\ \hline
 
  A_3+A_2+2A_1
 & \scriptstyle{\checkmark}
 & \scriptstyle{\checkmark}
 \\ \hline
 
 D_4+3A_1
 & \scriptstyle{\checkmark}
 & \scriptstyle{\checkmark}
 \\ \hline
 
  2A_3+A_1
 & \scriptstyle{\checkmark}
 & \scriptstyle{\checkmark}
 \\ \hline

 A_4+A_2+A_1
 & \scriptstyle{\checkmark}
 & \scriptstyle{\checkmark}
 \\ \hline
 
  A_5+2A_1
 & \scriptstyle{\checkmark}
 & \scriptstyle{\checkmark}
 \\ \hline

  D_5+2A_1
 & \scriptstyle{\checkmark}
 & \scriptstyle{\checkmark}
 \\ \hline

   A_4+A_3
 & \scriptstyle{\checkmark}
 & \scriptstyle{\checkmark}
 \\ \hline

  D_4+A_3
 & \scriptstyle{\checkmark}
 & \scriptstyle{\checkmark}
 \\ \hline
 
     \end{array}
 \hspace{1.5mm}
 \begin{array}{|c||c|c|} 
 \hline
 
\multirow{2}{*}{$\Gamma' \hookrightarrow E_8$} 
 & \multicolumn{2}{c|}{\text{occurs if}} \\
 & {\text{$p \neq 2$}}
 & {\text{$p = 2$}}
 \\ \hline \hline

 A_5+A_2
 & \scriptstyle{\checkmark}
 & \scriptstyle{\checkmark}
 \\ \hline
 
  D_5+A_2
 & \scriptstyle{\checkmark}
 & \scriptstyle{\checkmark}
 \\ \hline
 
  A_6+A_1
 & \scriptstyle{\checkmark}
 & \scriptstyle{\checkmark}
 \\ \hline
 
 D_6+A_1
 & \scriptstyle{\checkmark}
 & \scriptstyle{\checkmark}
 \\ \hline
 
 E_6+A_1
 & \scriptstyle{\checkmark}
 & \scriptstyle{\checkmark}
 \\ \hline
 
  A_7
 & \scriptstyle{\checkmark}
 & \scriptstyle{\checkmark}
 \\ \hline
 
  D_7
 & \scriptstyle{\checkmark}
 & \scriptstyle{\checkmark}
 \\ \hline
 
  E_7
 & \scriptstyle{\checkmark}
 & \scriptstyle{\checkmark}
 \\ \hline
 
  8A_1
 & \times
 & \scriptstyle{\checkmark}
 \\ \hline
 
 D_4+4A_1
 & \times
 & \scriptstyle{\checkmark}
 \\ \hline
 
  4A_2
 & \scriptstyle{\checkmark}
 & \scriptstyle{\checkmark}
 \\ \hline
 
  2A_3+2A_1
 & \scriptstyle{\checkmark}
 & \times
 \\ \hline
 
  A_5+A_2+A_1
 & \scriptstyle{\checkmark}
 & \scriptstyle{\checkmark}
 \\ \hline
 
  D_6+2A_1
 & \scriptstyle{\checkmark}
 & \scriptstyle{\checkmark}
 \\ \hline
 
    2A_4
 & \scriptstyle{\checkmark}
 & \scriptstyle{\checkmark}
 \\ \hline
 
   2D_4
 & \scriptstyle{\checkmark}
 & \scriptstyle{\checkmark}
 \\ \hline
 
   D_5+A_3
 & \scriptstyle{\checkmark}
 & \scriptstyle{\checkmark}
 \\ \hline
 
  E_6+A_2
 & \scriptstyle{\checkmark}
 & \scriptstyle{\checkmark}
 \\ \hline
 
   A_7+A_1
 & \scriptstyle{\checkmark}
 & \scriptstyle{\checkmark}
 \\ \hline
 
   E_7+A_1 
 & \scriptstyle{\checkmark}
 & \scriptstyle{\checkmark}
 \\ \hline
 
    A_8
 & \scriptstyle{\checkmark}
 & \scriptstyle{\checkmark}
 \\ \hline
 
   D_8
 & \scriptstyle{\checkmark}
 & \scriptstyle{\checkmark}
 \\ \hline
 
   E_8
 & \scriptstyle{\checkmark}
 & \scriptstyle{\checkmark}
 \\ \hline

 & 
 & 
 \\ \hline

 \end{array}$$
\vspace{-4mm}
\caption{$\Gamma' \subseteq E_8$ occurring on weak del Pezzo surfaces} 
\label{Table: Lattices embedding into E8}
\end{table}

\begin{Remark} \label{Remark exceptional cases KawaNaga}
We will see in Proposition \ref{Prop RDP iff RDP deg 1}, that, if $p \neq 2$, all rational double points that occur on some RDP del Pezzo surface also occur on an RDP del Pezzo surface of degree $1$. This fails for precisely four RDP configurations in characteristic $2$, namely for $\Gamma \in \{E_7^0,D_6^0 + A_1, D_4^0 + 3A_1, 7A_1\}$, each of which occurs on a unique RDP del Pezzo surface of degree $2$. These four surfaces coincide with the surfaces described in \cite[Theorem 1.4.(2)]{KawakamiNagaokaNeu-Pathologies}. Note, however, that the $(-2)$-curve configurations of types $E_7$ and $D_6 + A_1$ do occur on weak del Pezzo surfaces of degree $1$ in characteristic $2$, whereas the configurations of $(-2)$-curves $D_4 + 3A_1$ and $7A_1$ do not.
\end{Remark}

\begin{Remark}
One way of classifying RDP configurations on RDP del Pezzo surfaces over the complex numbers is to reduce to the case of RDP del Pezzo surfaces of Picard rank $1$ or $2$ as described in \cite[Lemma 2, Lemma 4]{MiyanishiZhang2}. However, as it is unclear whether this reduction also works in positive characteristic (in particular if $p = 2,3,5$), and since the classification of RDP del Pezzo surfaces of Picard rank $2$ in positive characteristic is not available yet, we will pursue a different approach in this article.
\end{Remark}

The structure of this article and thus also the structure of the proof of Theorem \ref{Thm main} is as follows: after recalling the classification of non-taut rational double points in Section \ref{Section Non-taut rational double points in positive characteristic}, we show in Section \ref{From del Pezzo surfaces to del Pezzo surfaces of degree $1$} that an RDP configuration occurs on an RDP del Pezzo surface if and only if it occurs on an RDP del Pezzo surface of degree $1$, with the exception of the four configurations mentioned in Remark \ref{Remark exceptional cases KawaNaga}. Then, in Section \ref{From del Pezzo surfaces of degree $1$ to rational (quasi-)elliptic surfaces}, we recall the well-known connection between RDP del Pezzo surfaces of degree $1$ and Weierstra{\ss} models of rational (quasi-)elliptic surfaces. In Section \ref{Configurations of $(-2)$-curves on weak del Pezzo surfaces: Reduction to non-taut RDPs}, we explain how this connection, and the theory of Mordell--Weil groups, can be exploited to classify all configurations of $(-2)$-curves that can occur on weak del Pezzo surfaces. This reduces Question \ref{Question central!} to non-taut RDPs and thus to characteristics $2,3,$ and $5$. Finally, the bulk of the article is devoted to the classification of RDP del Pezzo surfaces of degree $1$ with at least one non-taut rational double point in characteristic $2,3$, and $5$. This is achieved by using the classification of singular fibers of rational (quasi-)elliptic surfaces due to Ito \cite{QuasiEllipticChar3}, \cite{QuasiEllipticChar2} and Jarvis--Lang--Rimmasch--Rogers--Summers--Petrosyan \cite{EllipticChar3}, \cite{EllipticChar2} to derive simple equations for these RDP del Pezzo surfaces that allow us to explicitly determine the Artin coindices of the rational double points that occur.

\begin{Notationx}
In Tables \ref{TableD4D5-char2}, \ref{TableD6D7D8-char2}, \ref{TableE6E7E8-char2}, \ref{TableE6E7E8-char3}, and \ref{TableE8-char5} we list all possible RDP configurations containing a non-taut rational double point
in Column 1. In Column 2, we give simplified Weierstra{\ss} equations for all RDP del Pezzo surfaces of degree $1$ containing the respective configuration. If an extra condition on parameters in such an equation leads to extra RDPs, the condition is written under the respective equation separated by a dashed line. The equation above a dashed line is assumed to satisfy none of the conditions listed below it.
In Columns 3 and 4 we give the discriminant $\Delta$ and the $j$-invariant $j$ (see Subsection \ref{Subsec Tate algo} for explicit formulae) of the corresponding rational \mbox{(quasi-)}elliptic surface, whose type in the notation of Lang/Ito/Miranda--Persson is given in Column 5 and we note in Column 6 whether the fibration is elliptic or quasi-elliptic.
\end{Notationx}
 
\smallskip 
\subsection*{Acknowledgments}
The author would like to thank Gebhard Martin for interesting discussions, \mbox{Christian} Liedtke for comments on a first version of this article, as well as Igor Dolgachev and De-Qi Zhang for interesting remarks.
 
%\bigskip \bigskip

\section{Non-taut rational double points in positive characteristic} \label{Section Non-taut rational double points in positive characteristic}

Recall that \emph{rational double point} (RDP) is one of the names for a canonical surface singularity. The $n \geq 1$ irreducible components of the exceptional divisor of its minimal resolution span a negative definite root lattice $\Gamma_n$ of type $\Gamma \in \{A, D, E\}$.
If $p \neq 2,3,5$, every rational double point is \emph{taut}, that is, its formal isomorphism class is uniquely determined by $\Gamma_n$. It turns out that this fails for certain rational double points in small characteristics. Nevertheless, Artin \cite{ArtinCoveringsRDP} was able to give a classification of formal isomorphism classes of these non-taut rational double points and, in particular, he proved that there are only finitely many isomorphism classes $\Gamma_n^{1},\hdots,\Gamma_n^{k_n}$ for each fixed $\Gamma_n$. The number $k$ in $\Gamma_n^k$ is called \emph{Artin coindex} of the rational double point
(\textit{e.g.}, if $p = 5$, there are two distinct rational double points $E_8^0$ and $E_8^1$ both of which have resolution graph of type $E_8$).
We call a formal sum $\Gamma = \sum_i \Gamma_{i,n_i}^{k_i}$ of such RDPs $\Gamma_{i,n_i}^{k_i}$ an \emph{RDP configuration}.

In the following Table \ref{TableNonTaut-char235} (where $n \geq 2$ and $1 \leq r \leq n-1$), we listed Artin's equations for the non-taut rational double points together with the dimension $m$ of their miniversal deformation spaces. Here, we observe that the Artin coindices for a given Dynkin type can be distinguished by $m$:

\begin{Observation} \label{Obs RDPs distinguished by m}
The completions of two rational double points are isomorphic if and only if they have the same resolution graph and the dimensions $m$ of their miniversal deformation spaces coincide.
\end{Observation}

The following well-known Proposition \ref{Prop dimension(def)} provides a way to calculate $m$ for hypersurface singularities, and thus in particular for the rational double points.

\begin{Proposition} \cite[Chapter~1, \S 4]{ArtinDeftheory} \label{Prop dimension(def)}
Let the local ring $R= k[x,y,z]_{(x,y,z)}/(f)$ be a normal surface singularity given by one equation $f \in k[x,y,z]$. Then, the tangent space $T_f$ of the deformation functor ${\rm Def}_{\Spec R}$ of $\,{\Spec R}$ is given by 
\begin{equation*}
    T_f := {\rm Def}_{\Spec R} (k[\epsilon]/(\epsilon^2)) \cong \left.\raisebox{.3em}{$k[x,y,z]_{(x,y,z)}$} \middle/ \raisebox{-.3em}{$\left(f, \frac{\partial f}{\partial x}, \frac{\partial f}{\partial y}, \frac{\partial f}{\partial z}\right)$}\right. .
\end{equation*}
\end{Proposition}

\noindent
\begin{table}[h]
\renewcommand{\arraystretch}{1.1}
\noindent
\begin{minipage}{0.55\textwidth}
 $$
 \begin{array}{|c|c|c|}
  \multicolumn{3}{c}{\textbf{Characteristic 2}}\\
 \hline
\text{Type}
& {\text{Equation}}
& m
\\

\hline \hline

\multicolumn{3}{|c|}{\bf{D_n}}
\\
\hline

  D_{2n}^0  
    & z^2 + x^2y + xy^n
    & 4n
  \\
\hline

  D_{2n}^r 
  & z^2 + x^2y + xy^n + xy^{n-r}z
  & 4n-2r
  \\
\hline

  D_{2n+1}^0  
  & z^2 + x^2y + y^nz
  & 4n
  \\
\hline

  D_{2n+1}^r  
  & z^2 + x^2y + y^nz + xy^{n-r}z
& 4n-2r
  \\
\hline

\multicolumn{3}{|c|}{\bf{E_6}}
\\
\hline

  \text{\hspace{4mm}}E_6^0 \text{\hspace{4mm}}
  & {z^2 + x^3 + y^2z} 
  & \text{\hspace{5mm}}8\text{\hspace{5mm}}
  \\
\hline

  E_6^1
  & {z^2 + x^3 + y^2z + xyz}
  & 6
  \\
\hline

\multicolumn{3}{|c|}{\bf{E_7}}
\\
\hline

  E_7^0
  & {z^2 + x^3 + xy^3}
  & 14
  \\
\hline

  E_7^1
  & {z^2 + x^3 + xy^3 + x^2yz}
  & 12
  \\
\hline

  E_7^2
  & {z^2 + x^3 + xy^3 + y^3z}
  & 10
  \\
\hline

  E_7^3
  & {z^2 + x^3 + xy^3 + xyz}
  & 8
  \\
\hline

\multicolumn{3}{|c|}{\bf{E_8}}
\\
\hline

  E_8^0
  & \text{\hspace{12mm}} {z^2 + x^3 + y^5} \text{\hspace{12mm}}
  & 16
  \\
\hline

  E_8^1
  & {z^2 + x^3 + y^5 + xy^3z}
  & 14
  \\
\hline

  E_8^2
  & {z^2 + x^3 + y^5 + xy^2z}
  & 12
  \\
\hline

  E_8^3
  & {z^2 + x^3 + y^5 + y^3z}
  & 10
  \\
\hline

  E_8^4
  & {z^2 + x^3 + y^5 + xyz}
  & 8
  \\
\hline

\end{array}$$
\end{minipage}
\begin{minipage}{0.44\textwidth} 
 $$
 \begin{array}{|c|c|c|}
  \multicolumn{3}{c}{\textbf{Characteristic 3}}\\
 \hline
\text{Type}
& \text{Equation}
& m
\\

\hline \hline

\multicolumn{3}{|c|}{\bf{E_6}}
\\
\hline

  \text{\hspace{3mm}}E_6^0\text{\hspace{3mm}}
  & z^2 + x^3 + y^4
  & 9
  \\
\hline

  E_6^1
  & z^2 + x^3 + y^4 + x^2y^2 
  & 7
  \\
\hline

\multicolumn{3}{|c|}{\bf{E_7}}
\\
\hline

  E_7^0
  & z^2 + x^3 + xy^3
  & 9
  \\
\hline

  E_7^1
  & z^2 + x^3 + xy^3 + x^2y^2
  & 7
  \\
\hline

\multicolumn{3}{|c|}{\bf{E_8}}
\\
\hline

  E_8^0
  & \text{\hspace{10mm}} z^2 + x^3 + y^5 \text{\hspace{10mm}}
  & 12
  \\
\hline

  E_8^1
  & z^2 + x^3 + y^5 + x^2y^3
  & 10
  \\
\hline

  E_8^2
  & z^2 + x^3 + y^5 + x^2y^2
  & \text{\hspace{4mm}}8\text{\hspace{4mm}}
  \\
\hline

\end{array}$$
\bigskip
\smallskip
$$
 \begin{array}{|c|c|c|}
 \multicolumn{3}{c}{\textbf{Characteristic 5}}\\
 \hline
\text{Type}
& \text{Equation}
& m
\\

\hline \hline

\multicolumn{3}{|c|}{\bf{E_8}}
\\
\hline

  \text{\hspace{3mm}}E_8^0\text{\hspace{3mm}}
  & \text{\hspace{10mm}}z^2 + x^3 + y^5\text{\hspace{10mm}}
  & 10
  \\
\hline

  E_8^1
  & z^2 + x^3 + y^5 + xy^4
  & \text{\hspace{4mm}}8\text{\hspace{4mm}}
  \\
\hline

\end{array}$$
\end{minipage}
\caption{Types of non-taut rational double points in $\Char(k)=2,3,5$}
\label{TableNonTaut-char235}
\end{table}

\section{From del Pezzo surfaces to rational (quasi-)elliptic surfaces} \label{Section From del Pezzo surfaces to rational (quasi-)elliptic surfaces}

In this section, we reduce Question \ref{Question central!} to the corresponding question for Weierstra{\ss} models of rational (quasi-)elliptic surfaces.
On the way,  we recall the necessary background on del Pezzo surfaces and rational (quasi-)elliptic surfaces as well as their connection.

\subsection{From del Pezzo surfaces to del Pezzo surfaces of degree $1$} \label{From del Pezzo surfaces to del Pezzo surfaces of degree $1$}

Recall the following related notions of \emph{del Pezzo surfaces}.

\begin{Definition} \label{Def delPezzos}
Let $X$ and $\widetilde{X}$ be projective surfaces.
\begin{itemize}
\item
$X$ is a \emph{del Pezzo surface} if it is smooth and $-K_X$ is ample.
\item
$\widetilde{X}$ is a \emph{weak del Pezzo surface} if it is smooth and $-K_{\widetilde{X}}$ is big and nef.
The lattice $\sum_i \Gamma_{i,n_i} \subseteq \Pic (\widetilde{X})$ spanned by all the $(-2)$-curves on $\widetilde{X}$ is called \emph{configuration of} $(-2)$\emph{-curves on} $\widetilde{X}$.
\item 
$X$ is an \emph{RDP del Pezzo surface} if all its singularities are rational double points and $-K_X$ is ample.
The formal sum $\sum_i \Gamma_{i,n_i}^{k_i}$ of the formal isomorphism classes of all the rational double points on $X$ is called \emph{RDP configuration of} $X$.
\end{itemize}
In all the above cases, the number $\deg (X)= K_X^2$ (respectively $\deg (\widetilde{X})= K_{\widetilde{X}}^2$) is called the \emph{degree of} $X$ (respectively $\widetilde{X}$).
\end{Definition}

 Note that weak del Pezzo surfaces are precisely the minimal resolutions of RDP del Pezzo surfaces and every RDP del Pezzo surface $X$ is obtained by contracting all the $(-2)$-curves on a weak del Pezzo surface $\widetilde{X}$.
 If the RDP configuration of $X$ is $\sum_i \Gamma_{i,n_i}^{k_i}$, then the configuration of $(-2)$-curves on $\widetilde{X}$ is $\sum_{i} \Gamma_{i,n_i}$.

The following observation tells us that, in order to understand rational double points on RDP del Pezzo surfaces, it suffices to understand them on RDP del Pezzo surfaces of degree $1$ with precisely four exceptions in characteristic $2$:

\begin{Proposition} \label{Prop RDP iff RDP deg 1}
Let $\Gamma = \sum_i \Gamma_{i,n_i}^{k_i}$ be an RDP configuration.
If $\,\Gamma$ occurs on an RDP del Pezzo surface, but not on an RDP del Pezzo surface of degree $1$,
then $p = 2$ and $\Gamma$ is one of the following:
\begin{enumerate}
    \item[(A.)]\label{Exception A} $\Gamma = E_7^0$
    \item[(B.)]\label{Exception B} $\Gamma = D_6^0 + A_1$ 
    \item[(C.)]\label{Exception C} $\Gamma = D_4^0 + 3A_1$  
    \item[(D.)]\label{Exception D} $\Gamma = 7A_1$
\end{enumerate}
Moreover, there is a unique RDP del Pezzo surface (necessarily of degree $2$) realizing each of these exceptional cases.
\end{Proposition}

\begin{proof}
If $X$ is an RDP del Pezzo surface of degree $d \geq 2$ with RDP configuration $\Gamma$, we want to construct an RDP del Pezzo surface $X_1$ of degree $1$ with the same configuration $\Gamma$ by finding a point $p \in \widetilde{X}$, where $\widetilde{X}$ is the minimal resolution of $X$, such that ${\rm Bl}_p(\widetilde{X})$ is again a weak del Pezzo surface and such that ${\rm Bl}_p(\widetilde{X})$ and $\widetilde{X}$ have the same configuration of $(-2)$-curves. Contracting the $(-2)$-curves on ${\rm Bl}_p(\widetilde{X})$ yields an RDP del Pezzo surface of degree $(d-1)$ with RDP configuration $\Gamma$, so the claim will follow by induction.

We fix a realization of $\widetilde{X}$ as a blow-up of $\bbP^2$ in (possibly infinitely near) points
$p_1,\hdots,p_{9-d}$. Using the precise description of $(-2)$-curves on $\widetilde{X}$ (see \textit{e.g.}
\cite[Lemma $2.8.(i)$]{MartinStadlmayr}), we see that the blow-up of $\widetilde{X}$ in a point $p \in \widetilde{X}$
will be a weak del Pezzo surface with the same configuration of $(-2)$-curves as $\widetilde{X}$ if and only if the
following two conditions are satisfied:
\begin{enumerate}
    \item \label{Condition avoid negative curves} $p$ does not lie on a $(-1)$- or $(-2)$-curve on $\widetilde{X}$, and
    \item \label{Condition avoid singular cubics}if $d = 2$, then $p$ is not the singular point of the strict transform $\widetilde{C}$ of an irreducible singular cubic $C \subseteq \bbP^2$ through $p_1, \hdots, p_{7}$.
\end{enumerate}

\noindent Since there are only finitely many negative curves on $\widetilde{X}$, we can always find a $p$ that satisfies Condition (\ref{Condition avoid negative curves}). Thus, we may assume that $d = 2$.
To deal with Condition (\ref{Condition avoid singular cubics}), note that every $\widetilde{C}$ as in Condition (\ref{Condition avoid singular cubics}) is a member of the two-dimensional linear system $|-K_{\widetilde{X}}|$ and consider the variety 
\begin{equation*}
\mathfrak{I} = \{(\widetilde{C},p) \text{ }\mid\text{ } \widetilde{C}  \in |-K_{\widetilde{X}}| \text{ is integral and singular, and } p \text{ is its singular point}\} \text{.}
\end{equation*}

\noindent Denoting the sublocus of singular curves in $|-K_{\widetilde{X}}|$ as $|-K_{\widetilde{X}}|_{\rm sing}$, we have a correspondence
$$
\xymatrix{
\mathfrak{I} \ar[r]^{{\rm pr}_2} \ar[d]_{{\rm pr}_1} & \widetilde{X} \\
|-K_{\widetilde{X}}|_{\rm sing} \text{ } \ar@{^{(}->}[r] & |-K_{\widetilde{X}}|\text{,}
}
$$
\noindent 
where ${\rm pr}_1$ is quasi-finite. If the image of ${\rm pr}_2$ is not dense in $\widetilde{X}$, then we can find a $p$ satisfying Conditions (\ref{Condition avoid negative curves}) and (\ref{Condition avoid singular cubics}). Hence, we have to show that if ${\rm pr}_2$ is dominant, then $\Gamma$ is one of the four Exceptions (\hyperref[Exception A]{A.}), (\hyperref[Exception B]{B.}), (\hyperref[Exception C]{C.}), or (\hyperref[Exception D]{D.})

If ${\rm pr}_2$ is dominant, then 
$$
2 = \dim |-K_{\widetilde{X}}| \geq \dim |-K_{\widetilde{X}}|_{\rm sing} = \dim \mathfrak{I} \geq \dim (\widetilde{X})=2 \text{,}
$$
and since $|-K_{\widetilde{X}}|_{\rm sing}$ is closed in $|-K_{\widetilde{X}}| \cong \bbP^2$, we have $|-K_{\widetilde{X}}|_{\rm sing} = |-K_{\widetilde{X}}|$, that is, every anti-canonical curve on $\widetilde{X}$ is singular.
By \cite[Theorem $1.4.(2)$]{KawakamiNagaokaNeu-Pathologies} the $(-2)$-curve configuration $\Gamma'$ associated to $\Gamma$ is one of the following:
\begin{enumerate}
    \item[($\mathrm{A}'$.)]\label{Exception A'} $\Gamma' = E_7$ 
    \item[($\mathrm{B}'$.)]\label{Exception B'} $\Gamma' = D_6 + A_1$  
    \item[($\mathrm{C}'$.)]\label{Exception C'} $\Gamma' = D_4 + 3A_1$ 
    \item[($\mathrm{D}'$.)]\label{Exception D'} $\Gamma' = 7A_1$ 
\end{enumerate}
So, we still have to determine the Artin coindices in Cases (\hyperref[Exception A']{$\mathrm{A}'$.}), (\hyperref[Exception B']{$\mathrm{B}'$.}) and (\hyperref[Exception C']{$\mathrm{C}'$.}). By \cite[Remark 5.2, Remark 5.5, and Remark 5.8]{KawakamiNagaokaNeu-Pathologies} and \cite[Theorem 5.2.]{QuasiEllipticChar2}, the non-taut rational double points in these cases are isomorphic to RDPs that occur in partial resolutions of the affine surfaces in $\mathbb{A}^3$ given by
\begin{enumerate}
\item[($\mathrm{A}''.$)]\label{Exception A''}  $y^2 = x^3 + t^5$,
\item[($\mathrm{B}''.$)]\label{Exception B''}  $y^2 = x^3 + t^3x$, and
\item[($\mathrm{C}''.$)]\label{Exception C''}  $y^2 = x^3 + (t^3 + t)x$, respectively.
\end{enumerate}
 We see from Table \ref{TableNonTaut-char235} that the first two cases are $E_8^0$ and $E_7^0$, respectively. In Case (\hyperref[Exception C'']{$\mathrm{C}''.$}) one can apply Proposition \ref{Prop dimension(def)} to see that the singularity is of type $D_6^0$. By \cite[Theorem 2.70., Table 10]{Stadlmayr-Barbeit}, the Artin coindex of
a rational double point appearing in a partial resolution of a rational double point with coindex $0$ is itself $0$. Hence, the Artin coindices in the exceptional Cases (\hyperref[Exception A']{$\mathrm{A}'.$}), (\hyperref[Exception B']{$\mathrm{B}'.$}), and (\hyperref[Exception C']{$\mathrm{C}'.$}) are all $0$. 

 The existence and uniqueness of the four exceptional cases was proved in \cite[Theorem 1.4, Table~1]{KawakamiNagaokaNeu-Pathologies}.
\end{proof}
\medskip

\noindent Thus, Proposition \ref{Prop RDP iff RDP deg 1} reduces the initial Question \ref{Question central!} to the following one: 

\begin{Question} \label{Question deg 1}
Which rational double points occur on RDP del Pezzo surfaces of degree $1$?
\end{Question}

\subsection{From del Pezzo surfaces of degree $1$ to rational (quasi-)elliptic surfaces} \label{From del Pezzo surfaces of degree $1$ to rational (quasi-)elliptic surfaces}
To answer Question \ref{Question deg 1}, we will exploit the connection between weak (resp. RDP) del Pezzo surfaces of degree $1$ and rational \mbox{(quasi-)}elliptic surfaces (resp. Weierstra{\ss} models of those). For this, let us first recall their definition (see for example \cite[Chapter V]{DolgachevCossec}).

\begin{Definition}
Let $Y$ and $\widetilde{Y}$ be projective surfaces.
\begin{itemize}
    \item $\widetilde{Y}$ is a \emph{rational (quasi-)elliptic surface} if it is smooth, rational, and admits a morphism $f: \widetilde{Y} \to \bbP^1$ such that the following conditions hold: 
    \begin{itemize}
        \item $f$ is surjective with $f_*\cal{O}_{\widetilde{Y}} = \cal{O}_{\bbP^1}$,
        \item the generic fiber of $f$ is a regular curve of arithmetic genus $1$,
        \item there are no $(-1)$-curves in fibers of $f$, and
        \item $f$ admits a section $\sigma_0: \bbP^1 \to \widetilde{Y}$.
    \end{itemize}
Moreover, the group ${\rm MW}(f)$ of sections of $f$ is called \emph{Mordell--Weil group} of $f: \widetilde{Y} \to \bbP^1$.
    \item The \emph{Weierstra{\ss} model} $Y$ \emph{of} $\widetilde{Y}$ is the surface obtained from $\widetilde{Y}$ by contracting all components of fibers of $f$ that do not meet $\sigma_0(\bbP^1)$.
\end{itemize}
\end{Definition}

\begin{Remark}
In the literature one usually finds the definition of a (quasi-)elliptic surface as a pair of a surface and a (quasi-)elliptic fibration. Since $\widetilde{Y}$ is rational, the canonical bundle formula shows that it admits a unique (quasi-)elliptic fibration induced by $|-K_{\widetilde{Y}}|$, so we do not need to specify the fibration. Similarly, while a priori the Weierstra{\ss} model $Y$ seems to depend on the chosen section $\sigma_0$ of $f$, any two such sections are interchanged by an automorphism of $\widetilde{Y}$, so the associated Weierstra{\ss} models are isomorphic, and thus we will not keep track of the section.
\end{Remark}

Note that, because all components of reducible fibers of $f$ are $(-2)$-curves, the Weierstra{\ss} model $Y$ of $\widetilde{Y}$ has only rational double points as singularities. So, analogously to Section \ref{From del Pezzo surfaces to del Pezzo surfaces of degree $1$}, we define the \emph{configuration of} $(-2)$\emph{-curves on} $\widetilde{Y}$
and the \emph{RDP configuration of} $Y$. Note, however, that the configuration of $(-2)$-curves on $\widetilde{Y}$ is not a sum of root lattices in general.

We assume that the reader is familiar with the Kodaira--N\'eron classification of singular fibers of \mbox{(quasi-)} elliptic surfaces as described for example in \cite[Table 4.1, p.365]{Silverman2-AdvancedTopics}. In the following Table \ref{Table KodairaNeron and RDPs}, we summarize which Kodaira--N\'eron type of a fiber in $\widetilde{Y}$ leads to which rational double point on its image in $Y$. Here, we denote a smooth point by $A_0$, and we have $n \geq 1$ for type ${\rm I}_n$ and $n \geq 0$ for type ${\rm I}_n^*$.

\begin{table}[h]
\renewcommand{\arraystretch}{1.2}
    \centering
    $ \begin{array}{|c||c|c|c|c|c|c|c|c|c|}
    \hline
        \text{Kodaira--N\'eron type} &  {\rm I}_0 & {\rm I}_n  & {\rm II} & {\rm III} & {\rm IV}  & {\rm I}_n^* & {\rm IV}^* & {\rm III}^* & {\rm II}^* \\
        \hline \hline
       \text{Rational double point}  & A_0 & A_{n-1} & A_0 & A_1 & A_2  & D_{4+n} & E_6 & E_7 & E_8 \\
       \hline
    \end{array} $
    %\vspace{-3mm}
    \caption{Kodaira--N\'eron types and corresponding rational double points} 
    \label{Table KodairaNeron and RDPs}
\end{table}

By \cite[Section 8.3.2]{Dolgachev-classical} an RDP del Pezzo surface $X$ of degree $1$ is isomorphic to a sextic hypersurface $V(f_6) \subseteq \bbP(1,1,2,3)$ and, conversely, a sextic hypersurface in $\bbP(1,1,2,3)$ with at worst rational double point singularities defines an RDP del Pezzo surface of degree $1$. Such sextics are of the form 
\begin{equation} 
  \tag{W}  \label{W}
  y^2 \hspace{2mm}+ \hspace{2mm}a_1 xy\hspace{2mm} +\hspace{2mm} a_3 y \hspace{2mm}=\hspace{2mm} x^3\hspace{2mm} +\hspace{2mm} a_2 x^2 \hspace{2mm}+\hspace{2mm} a_4 x \hspace{2mm}+ \hspace{2mm}a_6,
\end{equation}
where the $a_i \in k[t,s]$ are homogeneous of degree $i$ and $t,s,x,y$ are of degrees $1,1,2,3$, respectively. 

Projecting $\bbP(1,1,2,3)$ onto $t$ and $s$ yields a rational map 
$\bbP(1,1,2,3) \dashrightarrow \bbP^1$, which, when restricted to the RDP del Pezzo $X= V(f_6)$ is given by the linear system $|-K_X|$ and has precisely one base point (at $s=t=0, y^2=x^3$). Blowing up the base point yields the Weierstra{\ss} model $Y \to \bbP^1$ of a rational \mbox{(quasi-)} elliptic surface $f: \widetilde{Y} \to \bbP^1$, where the zero section $\sigma_0$ on $Y$ resp. $\widetilde{Y}$ is the exceptional $(-1)$-curve of this blow-up. Conversely, for a rational (quasi-)elliptic surface $f: \widetilde{Y} \to \bbP^1$ with chosen section $\sigma_0$, contracting all components of fibers not meeting $\sigma_0$ yields its Weierstra{\ss} model $Y$, and contracting also $\sigma_0$, we obtain an RDP del Pezzo surface $X$ of degree $1$. In turn, $X$ is the anti-canonical model of a weak del Pezzo surface $\widetilde{X}$ of degree $1$, which, when blown-up in the base point of its anti-canonical linear system $|-K_{\widetilde{X}}|$ gives back $\widetilde{Y}$. This connection is summarized in the following commutative diagram:
\begin{equation}
 \label{Diagram connection}
\xymatrix{
  \widetilde{Y} \ar[d] \ar[r] \ar@/^1pc/[rr]^f & Y \ar[d] \ar[r] & \bbP^1  
\\
  \widetilde{X}  \ar[r] & X  \ar@{-->}[ur]&  
} 
\end{equation}

In particular, since the morphism $Y \to X$ is the blow-up of a smooth point, this diagram shows the following:

\begin{Observation} \label{Obs RDP deg 1 iff Wmodel}
A configuration of rational double points occurs on an RDP del Pezzo surface $X$ of degree $1$ if and only if it occurs on the Weierstra{\ss} model of a rational (quasi-)elliptic surface.
\end{Observation}

Thus, Question \ref{Question deg 1} is equivalent to the following one:

\begin{Question}
Which rational double points occur on Weierstra{\ss} models of rational (quasi-)elliptic surfaces?
\end{Question}

\subsection{Configurations of $(-2)$-curves on weak del Pezzo surfaces: reduction to non-taut RDPs} \label{Configurations of $(-2)$-curves on weak del Pezzo surfaces: Reduction to non-taut RDPs}
Let $\Gamma= \sum_i \Gamma_{i,n_i}$ be a configuration of $(-2)$-curves. Consider the following three conditions, where $p = \Char(k)$ and $\ell \neq p$ is a prime and the respective torsion parts are denoted by $[\cdot ]$.

\begin{equation*}
    \tag{E8}\label{E8}
    \text{There is an embedding } \iota: \Gamma \hookrightarrow E_8\text{.}
\end{equation*} 
\begin{equation*}
    \tag{E8+T[$\ell$]}\label{E8+T[l]}
    \text{(\ref{E8}) and } (E_8/ \iota(\Gamma))[\ell] \subseteq (\bbZ/ \ell \bbZ)^2 \text{.}
\end{equation*} 
\begin{equation*}
    \tag{E8+T[p]}\label{E8+T[p]}
    \text{(\ref{E8}), and if $p > 0$, } (E_8/ \iota(\Gamma))[p] \subseteq \bbZ/ p \bbZ \text{ or } (E_8/ \iota(\Gamma)) \cong (\bbZ/p \bbZ)^n \text{ for some } n \geq 0\text{.}
\end{equation*}

The reason why we consider the above three conditions is the following lemma.

\begin{Lemma} \label{Lemma if gamma occurs, then conditions satisfied}
Let $\Gamma= \sum_i \Gamma_{i,n_i}$ be a configuration of $(-2)$-curves. If $\Gamma$ occurs on a weak del Pezzo surface of degree $1$, then it satisfies Condition \emph{(\ref{E8+T[p]})} and \emph{(\ref{E8+T[l]})}
for all $\ell \neq p = \Char(k)$.
\end{Lemma}

\begin{proof}
Let $\widetilde{X}$ be a weak del Pezzo surface of degree $1$ realizing $\Gamma$. Let $f: \widetilde{Y} \to \bbP^1$ be the associated rational (quasi-)elliptic surface as in Diagram (\ref{Diagram connection}) with Mordell--Weil group ${\rm MW}(f)$.

By \cite[Theorem 3.1]{OguisoShioda}, there is an embedding $\iota: \Gamma \hookrightarrow E_8$ such that
\[{\rm rk}({\rm MW}(f)) = 8 - {\rm rk}(\Gamma) = {\rm rk} (E_8/\iota(\Gamma))\quad \text{and}\quad {\rm MW}(f)_{{\rm tors}} = (E_8/\iota(\Gamma))_{{\rm tors}}.\]
For all $n \geq 0$, we have 
$(E_8/\iota(\Gamma))[n] = ({\rm MW}(f))[n] \subseteq X_{\bar{\eta}}^{sm}[n]$, 
where $X_{\bar{\eta}}^{sm}$ is the smooth locus of the geometric generic fiber of $f$. If $f$ is elliptic, then $X_{\bar{\eta}}^{sm}[\ell] \cong (\bbZ/ \ell \bbZ)^2$ and $X_{\bar{\eta}}^{sm}[p] \subseteq \bbZ/ p\bbZ$, so the Conditions {(\ref{E8+T[p]})}, and {(\ref{E8+T[l]})} are satisfied for all $\ell \neq p = \Char(k)$. If $f$ is quasi-elliptic, then ${\rm MW}(f)$ is a finitely generated subgroup of $X_{\bar{\eta}}^{sm} \cong \mathbb{G}_a$, so ${\rm MW}(f) \cong (\bbZ/p\bbZ)^n$ for some $n \geq 0$. In particular, $E_8/\iota(\Gamma) \cong (\bbZ/p\bbZ)^n$, so again both conditions are satisfied.
\end{proof}
\medskip

The root sublattices of $E_8$ have been classified by Dynkin \cite[\S 5., Table 11, p.385]{DynkinSemisimple}. We can easily check which of them satisfy the conditions above.

\begin{Lemma} \label{Lemma which of all the configurations satisfy which conditions}
Let $\Gamma= \sum_i \Gamma_{i,n_i} \subseteq E_8$ be a configuration of $(-2)$-curves. Then, the following hold:
\begin{enumerate}
    \item $\Gamma$ satisfies \emph{({E8+T[q]})} for $q \neq 2$,
    \item If $p \neq 2$, then \emph{({E8+T[$\ell$ = 2]})} is satisfied if and only if $\Gamma \not \in \{D_4 + 4A_1, 8A_1, 7A_1\}$.
    \item If $p = 2$, then \emph{({E8+T[p=2]})} is satisfied if and only if $\Gamma \not \in \{D_4 + 3A_1, 2A_3 + 2A_1,  A_3 + 4A_1, 7A_1,6A_1\}$.
\end{enumerate}
\end{Lemma}

\prf
The groups $(E/\Gamma)_{\rm tors}$ have been calculated by Oguiso and Shioda in \cite[Corollary 2.1]{OguisoShioda} for all $\Gamma$ except $D_4 + 4A_1, 8A_1,$ and $7A_1$. Using their results, we leave it to the reader to check that for $\Gamma \not \in  \{D_4 + 4A_1, 8A_1, 7A_1\}$
the Conditions {(\ref{E8+T[l]})} are satisfied for all $\ell$, and $\Gamma$ does not satisfy Condition {(\ref{E8+T[p]})} if and only if $p = 2$ and $\Gamma \in \{D_4 + 3A_1, A_3 + 4A_1, 2A_3 + 2A_1, 6A_1\}$.
We will now treat the three remaining cases.
\begin{itemize}
    \item If $\Gamma = D_4 + 4A_1$, then $E_8/\Gamma = (\bbZ/2\bbZ)^3$ by \cite[Table 1]{QuasiEllipticChar2}. So, $\Gamma$ satisfies {(\ref{E8+T[p]})}, and it satisfies 
{(\ref{E8+T[l]})} if and only if $\ell \neq 2$.

\item If $\Gamma = 8A_1$, then $E_8/\Gamma = (\bbZ/2\bbZ)^4$ by \cite[Table 1]{QuasiEllipticChar2}.  So, $\Gamma$ satisfies {(\ref{E8+T[p]})}, and it satisfies 
{(\ref{E8+T[l]})} if and only if $\ell \neq 2$.

\item Finally, if $\Gamma = 7A_1$, then there is a unique embedding of $7A_1$ into $E_8$ by \cite[\S 5., Table 11, p.385]{DynkinSemisimple}. So, this embedding coincides with $7A_1 \hookrightarrow 8A_1 \hookrightarrow E_8$ and thus $E_8/7A_1$ has rank $1$ and $E_8/7A_1[2]$ contains $(\bbZ/2\bbZ)^3$. Hence, $\Gamma$ never satisfies ({E8+T[2]}), but it satisfies all \mbox{({E8 + T[q]})} with $q \neq 2$.
\end{itemize}
\vspace{-2mm}
\qed

\begin{Proposition} \label{Prop configurations of (-2)}
Let $\Gamma:= \sum_{i} \Gamma_{i,n_i}$ be a configuration of $(-2)$-curves. Then, the following are equivalent:
\begin{enumerate}
    \item \label{labelThm -2 confi occurs} $\Gamma$ occurs on a weak del Pezzo surface of degree $1$.
    \item \label{labelThm -2 confi satisfies condi}$\Gamma$ satisfies \emph{({E8+T[2]})}.
\end{enumerate}
\end{Proposition}

\prf
The implication (\ref{labelThm -2 confi occurs}) to (\ref{labelThm -2 confi satisfies condi}) follows immediately from Lemma \ref{Lemma if gamma occurs, then conditions satisfied}. 

For the converse, we have to show that every configuration of $(-2)$-curves that satisfies {({E8+T[2]})} occurs on a weak del Pezzo surface of degree $1$, or, equivalently, as the configuration of $(-2)$-curves that do not meet a fixed chosen section $\sigma_0$ on a rational (quasi-)elliptic surface. We remind the reader that we summarized the relation between these curve configurations and the corresponding Kodaira--N\'eron types in Table \ref{Table KodairaNeron and RDPs}.

For $p \neq 2,3$, this is precisely the content of \cite[Theorem 8.9]{SchuettShioda-MWlattices} (see also \cite[Remark 2.7]{OguisoShioda}).
For $p = 3$, it follows from the classification of singular fibers of elliptic surfaces in characteristic $3$ \cite{EllipticChar3}, that every $\Gamma$ that satisfies ({E8+T[2]}) occurs on an elliptic surface, except $\Gamma = 4A_2$. By \cite[Theorem 3.3]{QuasiEllipticChar3}, this $\Gamma$ is realized on a quasi-elliptic surface. Similarly, if $p = 2$, one can use \cite{EllipticChar2} and \cite{QuasiEllipticChar2} to check that every $\Gamma$ that satisfies ({E8+T[2]}) occurs on some (quasi-)elliptic surface.
\qed
\medskip

Combining the results of this section, we can give a proof of Theorem \ref{Thm main} for configurations of taut rational double points. In particular, this proves Theorem \ref{Thm main} in characteristic different from $2,3$, and $5$:

\medskip
\begin{proof}[Proof of Theorem \ref{Thm main} for taut RDPs]
Let $\Gamma = \sum_i \Gamma_{i,n_i}^{k_i}$ be an RDP configuration and assume that all the $\Gamma_{i,n_i}^{k_i}$ are taut, so $\Gamma$ is uniquely determined by its associated configuration of $(-2)$-curves $\Gamma' := \sum_{i} \Gamma_{i,n_i}$. Further assume that $\Gamma \neq 7A_1$ if $p=2$. Then, we have the following equivalences:

\noindent
\begin{center}
\begin{tabular}{p{4cm} c  l}
$\Gamma$ occurs on an RDP del Pezzo surface& $\xLeftrightarrow{\ \text{Prop. }\ref{Prop RDP iff RDP deg 1}\ }$ 
& $\Gamma$ occurs on an RDP del Pezzo surface of degree $1$\\
& $\xLeftrightarrow{\text{all } \Gamma_{i,n_i}^{k_i} \text{ taut}}$ & $\Gamma'$ occurs on a weak del Pezzo surface of degree $1$\\
& $\xLeftrightarrow{\ \text{Prop. }\ref{Prop configurations of (-2)}\, }$ & $\Gamma'$ satisfies ({E8+T[2]}) \\
& $\xLeftrightarrow{\ \, \text{Lem. }\ref{Lemma which of all the configurations satisfy which conditions}\ \, }$ & $\Gamma'$ satisfies (\hyperref[E8]{E8}) and \\
& &
$\Gamma' \not \in  \begin{cases}
\{D_4 + 4A_1,8A_1, 7A_1\} & \text{if } p \neq 2 \text{,} \\
\{2A_3 + 2A_1, A_3 + 4A_1, 6A_1\} & \text{if }p = 2 \text{.}
\end{cases}$
\end{tabular}
\end{center}

\noindent
Note that for the last equivalence we did not have to consider $D_4+3A_1$ since $D_4$ is not taut if $p=2$.\\
\indent
Since we already know by Proposition \ref{Prop RDP iff RDP deg 1} that $\Gamma= 7A_1$ occurs on an RDP del Pezzo surface of degree $2$, this proves Theorem \ref{Thm main} in the case where all the $\Gamma_{i,n_i}^{k_i}$ are taut.
\end{proof}
\smallskip

Thus, we have reduced Question \ref{Question deg 1} to the following:

\begin{Question}
Which non-taut RDPs occur on Weierstra{\ss} models of rational (quasi-)elliptic surfaces?
\end{Question}

The remainder of this article will be devoted to finding an answer to this question. 
 
 \section{Classification of non-taut RDP del Pezzo surfaces of degree $1$} \label{Classification of non-taut RDP del Pezzo surfaces of degree $1$}
 In this section, we will explain how to derive simple equations for non-taut RDP del Pezzo surfaces of degree $1$ from the classification of singular fibers of rational (quasi-)elliptic surfaces given in \cite{EllipticChar0}, \cite{EllipticChar3}, \cite{EllipticChar2}, \cite{QuasiEllipticChar3}, and \cite{QuasiEllipticChar2}. In the (quasi-)elliptic case, the equations found by Ito are already simplified, which is why we will focus on the elliptic case in Subsections \ref{Subsec Tate algo}, \ref{Simplified Weierstrass equations} and \ref{Subsec Autos preserving simplified Weierstrass equations}.
Using these simplified equations, it is straightforward to determine the Artin coindices of the rational double points that occur by applying Proposition \ref{Prop dimension(def)}. We will exemplify this final step in Example \ref{Ex dim def char 3} and Example \ref{Ex dim def char 2} in characteristic $2$ and $3$, respectively.

\subsection{Tate's algorithm for determining the type of a singular fiber in an elliptic pencil} \label{Subsec Tate algo}

Given a Weierstra{\ss} equation for the Weierstra{\ss} model of a rational elliptic surface, it is fairly standard to determine the Kodaira--N{\'e}ron types (see Table \ref{Table KodairaNeron and RDPs}) of its singular fibers by carrying out Tate's algorithm \cite{TateAlgo}. For the sake of self-containedness, we quickly recall Tate's algorithm: Let
\begin{equation} 
  \tag{W} 
  y^2 \hspace{2mm}+ \hspace{2mm}a_1 xy\hspace{2mm} +\hspace{2mm} a_3 y \hspace{2mm}=\hspace{2mm} x^3\hspace{2mm} +\hspace{2mm} a_2 x^2 \hspace{2mm}+\hspace{2mm} a_4 x \hspace{2mm}+ \hspace{2mm}a_6,
\end{equation}
be a Weierstra{\ss} equation of an elliptic curve over the function field $k(C)$ of a smooth curve $C$. Choose a closed point $c \in C$ and write $\widehat{\mathcal{O}}_{C,c} \cong k [[ t ]]$, assume $a_i \in k[[t]]$ and write $\nu(\cdot)={\rm ord}_t(\cdot)$ for the valuation on $k[[t]]$. Then, Tate defines the following quantities: \hfill
$$
\begin{array}{l}
b_2  :=  a_1^2 + 4a_2, \hspace{3mm}
b_4  :=   a_1a_3 + 2a_4, \hspace{3mm}
b_6  :=  a_3^2 + 4a_6 , \hspace{3mm}
b_8  :=   a_1^2a_6 - a_1a_3a_4 + 4a_2a_6 + a_2a_3^2 - a_4^2, \\
c_4  :=  b_2^2 - 24b_4,  \hspace{3mm}
c_6  :=  -b_2^3 + 36b_2b_4 - 216b_6 ,  \hspace{3mm}
\Delta  =   -b_2^2b_8 - 8b_4^3 - 27b_6^2 + 9b_2b_4b_6 \neq 0,  \hspace{3mm}
j  =   \frac{c_4^3}{\Delta}.
\end{array}
$$

Then, excluding the subalgorithm for determining $v >0$ in Step $7$ if $p=2$ (see \cite[p.50-51]{TateAlgo} for details in this case), Tate's algorithm is as follows, where $F$ denotes the fiber over $t=0$.

\begin{Algorithmx} \hfill
\begin{enumerate}[leftmargin=2cm]
\item [Step $1$:]
If $t \nmid \Delta$, then $F$ is of type ${\rm I}_0$. Else$\hdots$
\item [Step $2$:]
Change coordinates such that $t \mid a_3, t\mid a_4$ and $t\mid a_6$. \\ 
\noindent
If $t \nmid b_2$, then $F$ is of type ${\rm I}_v$ for $v= \nu(\Delta)$. Else$\hdots$
\item [Step $3$:]
If $t^2 \nmid a_6$, then $F$ is of type ${\rm II}$. Else$\hdots$
\item [Step $4$:]
If $t^3 \nmid b_8 $, then $F$ is of type ${\rm III}$. Else$\hdots$
\item [Step $5$:]
If $t^3 \nmid b_6 $, then $F$ is of type ${\rm IV}$. Else$\hdots$
\item [Step $6$:]
Change coordinates such that $t \mid a_1, t \mid a_2, t^2 \mid a_3, t^2 \mid a_4$ and $ t^3 \mid a_6$. \\
Consider the polynomial $P(T)=T^3 + a_2 \frac{T^2}{t} + a_4 \frac{T}{t^2} + a_6 \frac{1}{t^3}$. \\
\noindent
If $P$ has three distinct roots, then $F$ is of type ${\rm I}_0^*$. Else$\hdots$
\item [Step $7$:]
If $P$ has one single or one double root, then $F$ is of type $I_v^*$ for some $v >0$ (if $p \neq 2$, then $v= \nu(\Delta) -6$). Else$\hdots$
\item [Step $8$:]
Change variables such that the triple root is $0$, and $t^2 \mid a_2, t^3 \mid a_4$ and $t^4 \mid a_6$.  \\
Consider the polynomial $Q(Y)= Y^2 + \frac{a_3}{t^2} Y - a_{6} \frac{1}{t^4}$. \\
If $Q$ has distinct roots, then $F$ is of type ${\rm IV}^*$. Else$\hdots$
\item [Step $9$:]
Change variables such that the double root is $0$ and $t^3 \mid a_3$ and $ t^5 \mid a_6$.\\ 
If $t^4\nmid a_4 $, then $F$ is of type ${\rm III}^*$. Else$\hdots$
\item [Step $10$:]
If $t^6\nmid a_6 $, then $F$ is of type ${\rm II}^*$. Else$\hdots$
\item [Step $11$:]
Divide each $a_i$ by $t^i$ and repeat from Step $1$.
\end{enumerate}
\end{Algorithmx}

\subsection{Simplified Weierstra{\ss} equations} \label{Simplified Weierstrass equations}

Depending on the characteristic, the Weierstra{\ss} equation (\hyperref[W]{W}) (if $p=2$, under the additional assumption that it contains a non-taut rational double point) for an RDP del Pezzo surface $X$ of degree $1$ with associated rational elliptic surface $f: \widetilde{Y} \to \bbP^1$ can be simplified to an equation of the following form: 

\vspace{2mm}
\noindent
$\begin{array}{l c c c c c c c c c c c c c c l}
     \text{(W0)\hspace{0.1cm}}\label{W0} && y^2 &&& &&= &x^3 && &+& a_4 x &+& a_6 & \text{if } p \neq 2,3,  \\
     \text{(W3)\hspace{0.1cm}}\label{W3}& & y^2 &&& &&= &x^3 &+& a_2 x^2 &+& a_4 x &+& a_6 & \text{if } p =3,\text{ and} \\
     \text{(W2)\hspace{0.1cm}}\label{W2} &\multirow{2}{*}{\text{or}}& y^2 &+& a_1 xy && &= &x^3 & +&a_2 x^2&+& a_4 x &+& a_6 &  \multirow{2}{*}{\text{if } $p=2$\text{.}}\\
     \text{(W2$'$)\hspace{0.1cm}}\label{W2'} & & y^2 &&  &+& a_3 y&= &x^3 & +&a_2 x^2&+& a_4 x &+& a_6 &  \\
\end{array}$

\vskip\baselineskip
This is well-known if $p \neq 2$. If $p = 2$, to see that we can simplify (\hyperref[W]{W}) to an equation of the form
(\hyperref[W2]{W2}) or (\hyperref[W2']{W2$'$}), we first observe that by Table \ref{Table KodairaNeron and RDPs} the non-taut rational double points (see Table \ref{TableNonTaut-char235}) correspond to certain additive fibers of $f$. Hence, Tate's algorithm shows that, if $X$ contains a non-taut rational double point, then in (\hyperref[W]{W}) we may assume $t \mid a_1$ and $t \mid a_3$. Thus, if $a_1=0$, we get (\hyperref[W2']{W2$'$}) and if $a_1 \neq 0$, we can assume that $\frac{a_3}{a_1}$ is a polynomial and $x \mapsto x + \frac{a_3}{a_1}$ transforms (\hyperref[W]{W}) to an equation of the form (\hyperref[W2]{W2}).

\subsection{Automorphisms of $\bbP(1,1,2,3)$ preserving simplified Weierstra{\ss} equations} \label{Subsec Autos preserving simplified Weierstrass equations}

In order to find simple equations for degree $1$ RDP del Pezzo surfaces with non-taut rational double points, let us have a look at which automorphisms of $\bbP(1,1,2,3)$ send an equation of the form (\hyperref[W0]{W0}), (\hyperref[W3]{W3}), (\hyperref[W2]{W2}) or (\hyperref[W2']{W2$'$}) to an equation of the same form.

First, observe that substitutions in $t$ and $s$ only, always preserve these types of equations. Thus, let us focus on those automorphisms of $\bbP(1,1,2,3)$ fixing $t$ and $s$, that is, those inducing the trivial automorphism on $\bbP^1$.

\subsubsection{Automorphisms of $\bbP(1,1,2,3)$ over $\bbP^1$ if $p \neq 2,3$}

A general substitution fixing $t$ and $s$ and sending a Weierstra{\ss} equation of the form (\hyperref[W0]{W0}) to one of the same form is given by 
\begin{equation*}
x  \mapsto  \lambda^2 x, \hspace{3mm} y  \mapsto  \lambda^3 y \text{\hspace{4mm} with } \lambda \in k^* \text{.}
\end{equation*}
This sends (\hyperref[W0]{W0}) to $y^2 = x^3 + \frac{1}{\lambda^4}a_4x + \frac{1}{\lambda^6}a_6$.

\subsubsection{Automorphisms of $\bbP(1,1,2,3)$ over $\bbP^1$ if $p=3$} \label{Subsub general Auto-char3}

A general substitution fixing $t$ and $s$ and sending a Weierstra{\ss} equation of the form (\hyperref[W3]{W3}) to one of the same form is given by 
\begin{eqnarray*}
x & \mapsto & \lambda^2 x +f\\
y & \mapsto & \lambda^3 y 
\end{eqnarray*}
with $\lambda \in k^*$ and $f \in k[t,s]$ homogeneous of degree $2$. This sends (\hyperref[W3]{W3}) to 
\begin{equation*}
y^2 =x^3 + \frac{1}{\lambda^2} a_2 x^2 + \frac{1}{\lambda^4}(a_4+ 2a_2f)x + \frac{1}{\lambda^6}(a_6+ a_4f +a_2f^2 + f^3) \text{.}
\end{equation*}

\subsubsection{Automorphisms of $\bbP(1,1,2,3)$ over $\bbP^1$ if $p = 2$} \label{Subsub general Auto-char2}

\begin{enumerate}
    \item[(W2)]
    A general substitution fixing $t$ and $s$ and sending a Weierstra{\ss} equation of the form (\hyperref[W2]{W2}) to one of the same form is given by 
\begin{eqnarray*}
x & \mapsto & \lambda^2 x \\
y & \mapsto & \lambda^3 y + fx+ g
\end{eqnarray*}
with $\lambda \in k^*$ and $f,g \in k[t,s]$ homogeneous of degree $1$ and $3$, respectively. This sends (\hyperref[W2]{W2}) to 
\begin{equation*}
y^2 +  \frac{1}{\lambda}a_1 xy = x^3 + \frac{1}{\lambda^6}(\lambda^4 a_2 + \lambda^2 a_1 f + f^2)x^2 + \frac{1}{\lambda^4}(a_4+a_1g)x+ \frac{1}{\lambda^6}(a_6+g^2) \text{.}
\end{equation*}

    \item[(W2$'$)]
    A general substitution fixing $t$ and $s$ and sending a Weierstra{\ss} equation of the form (\hyperref[W2']{W2$'$}) to one of the same form is given by 
\begin{eqnarray*}
x & \mapsto & \lambda^2 x + f \\
y & \mapsto & \lambda^3 y + gx+ h
\end{eqnarray*}
with $\lambda \in k^*$ and $f,g,h \in k[t,s]$ homogeneous of degree $2,1$ and $3$, respectively. This sends (\hyperref[W2']{W2$'$}) to
\begin{equation*}
y^2 +  \frac{1}{\lambda^3}a_3y= x^3 + \frac{1}{\lambda^6}(\lambda^4 a_2 + g^2 + \lambda^4 f)x^2 + \frac{1}{\lambda^6}(\lambda^2 a_4 + a_3 g + \lambda^2 f^2)x + \frac{1}{\lambda^6}(a_6 + a_4f + a_3h + a_2 f^2 + f^3 + h^2) \text{.}
\end{equation*}
\end{enumerate}

\subsection{Proof of Theorem \ref{Thm main} in characteristic $5$}
Assume $p = 5$. Let $\Gamma = \sum_i \Gamma_{i,n_i}^{k_i}$ be an RDP configuration containing a non-taut RDP and let $\Gamma' = \sum_{i} \Gamma_{i,n_i}$ be the associated configuration of $(-2)$-curves. By Table \ref{TableNonTaut-char235}, we have $\Gamma_{i,n_i}= E_8$ for some $i$ and by Proposition \ref{Prop RDP iff RDP deg 1} and Lemma \ref{Lemma if gamma occurs, then conditions satisfied}, $\Gamma$ can only occur on an RDP del Pezzo surface if $\Gamma'$ embeds into $E_8$. Thus, to prove Theorem \ref{Thm main} in characteristic $5$, it suffices to consider $\Gamma \in \{E_8^0,E_8^1\}$. Note that $\Gamma'$ embeds into $E_8$ and $\Gamma' = E_8 \not \in \{D_4 + 4A_1,8A_1,7A_1\}$. On the other hand, the following proposition shows that both of these rational double points occur, so Theorem \ref{Thm main} holds in characteristic $5$.

\begin{Proposition} \label{Prop Equations-char5}
Each of the rational double points $E_8^0$ and $E_8^1$ occurs on an RDP del Pezzo surface $X$ of degree $1$. Moreover, every RDP del Pezzo surface of degree $1$ containing a non-taut rational double point is given by an equation as in Table \ref{TableE8-char5}.
\end{Proposition}

\prf
By Table \ref{Table KodairaNeron and RDPs}, the rational elliptic surface associated to an RDP del Pezzo surface of degree $1$ with a singularity of type $E_8$ admits a fiber of type ${\rm II}^*$. By \cite[Theorem 4.1.]{ExtremalCharpII} and \cite[Theorem~4.1., Tables~5.1 and~5.2]{ExtremalChar0} there
are precisely two such elliptic surfaces and their Weierstra{\ss} equations in $\bbP(1,1,2,3)$ are
\begin{equation}\label{eq char 5 0}
    y^2= x^3 + t^5s 
\end{equation}
and
\begin{equation}\label{eq char 5 1}
    y^2= x^3 + t^4x + t^5s \text{.}
\end{equation}
Considering the affine chart $s = 1$ and comparing with Table \ref{TableNonTaut-char235}, we see that Equation \eqref{eq char 5 0} defines a singularity of type $E_8^0$ and Equation \eqref{eq char 5 1} defines a singularity of type $E_8^1$.
\qed

\vspace{1mm}

\subsection{Proof of Theorem \ref{Thm main} in characteristic $3$} \label{Subsec Proof of Thm main char 3}
Assume $p = 3$. Let $\Gamma = \sum_i \Gamma_{i,n_i}^{k_i}$ be an RDP configuration containing a non-taut RDP and let $\Gamma' = \sum_{i} \Gamma_{i,n_i}$ be the associated configuration of $(-2)$-curves. By Table \ref{TableNonTaut-char235}, we have $\Gamma_{i,n_i} \in \{E_6,E_7,E_8\}$ for some $i$ and by Proposition \ref{Prop RDP iff RDP deg 1} and Lemma \ref{Lemma if gamma occurs, then conditions satisfied}, $\Gamma$ can only occur on an RDP del Pezzo surface if $\Gamma'$ embeds into $E_8$. Thus, by Dynkin's classification \cite[Table 11]{DynkinSemisimple}, to prove Theorem \ref{Thm main} in characteristic $3$, it suffices to consider 
\begin{equation} \label{eq possible types for Gamma in char 3}
\Gamma \in \{E_8^0,E_8^1, E_8^2, E_7^0 + A_1, E_7^0, E_7^1 + A_1, E_7^1, E_6^0 + A_2, E_6^0 + A_1, E_6^0, E_6^1 + A_2, E_6^1 + A_1, E_6^1\}\text{.}
\end{equation}

\noindent
Note that for all of the $\Gamma$ above, we have $\Gamma' \not \in \{D_4 +4A_1,8A_1,7A_1\}$. On the other hand, the following proposition shows that all these possiblities occur on some RDP del Pezzo surface, so Theorem \ref{Thm main} holds in \mbox{characteristic $3$.}

\begin{Proposition} \label{Prop Equations-char3}
Each RDP configuration $\Gamma$ in the List \eqref{eq possible types for Gamma in char 3} occurs on an RDP del Pezzo surface of degree $1$.
Moreover, every RDP del Pezzo surface of degree $1$ containing a non-taut rational double point admits an equation as in Table \ref{TableE6E7E8-char3}.
\end{Proposition}

\begin{proof}
By Table \ref{Table KodairaNeron and RDPs}, we have to study those RDP del Pezzo surfaces $X$ whose associated rational \mbox{(quasi-)} elliptic surface $\widetilde{Y}$ has a singular fiber of type ${\rm IV}^*$, ${\rm III}^*$, or ${\rm II }^*$.
In the elliptic case and in the notation of \cite{EllipticChar3}, these correspond to the Types $6A$, $6B$, $6C$, $7$, $8A$, and $8B$. In the quasi-elliptic case and in the notation of \cite[Theorem 3.3] {QuasiEllipticChar3}, these correspond to Cases (1) and (2).

All the Weierstra{\ss} equations for $X \subseteq \bbP(1,1,2,3)$ given in \cite{EllipticChar3} and \cite{QuasiEllipticChar3} are already of the form (\hyperref[W3]{W3}). To simplify them and determine the rational double points that occur, we will proceed along the following steps:
\begin{enumerate}
    \item[(1.)] Carry out a substitution in $t$ and $s$ only.
    \item[(2.)] Apply an automorphism of $\bbP(1,1,2,3)$ over $\bbP^1$ preserving the form (\hyperref[W3]{W3}) as in Subsection \ref{Subsub general Auto-char3}.
    \item[(3.)] Check for additional rational double points (\textit{e.g.} using Tate's algorithm (see Subsection \ref{Subsec Tate algo}) to determine the other reducible fibers of the underlying rational (quasi-)elliptic surface).
    \item[[(4.)] Determine the Artin coindices as described in Section \ref{Section Non-taut rational double points in positive characteristic}, \textit{e.g.} via Proposition \ref{Prop dimension(def)}.
 This will be left to the reader, but we will show how it works in Example \ref{Ex dim def char 3}.]

\end{enumerate}

\subsubsection*{Lang's Type 6A (${\rm IV}^*$)}
$X$ is given by $y^2= x^3+ c_0t^2x^2  + c_1t^3x + c_2t^4$ with $t \nmid c_1, t \nmid c_2$ and $c_i \in k[t,s]$ homogeneous of degree $i$. From now on, let us distinguish the cases $c_0=0$ and $c_0 \neq 0$.

\begin{itemize}
    \item[---] $c_0 = 0$:
    \begin{enumerate}
    \item[(1.)] Since $t \nmid c_1$, we can apply an automorphism of $\bbP^1$ to assume $c_1=s$. Then, scaling $s \mapsto \lambda^3 s, t \mapsto \lambda^{-1}t$ for an appropriate $\lambda$, we can write $c_2= s^2 + c_{2,1}ts + c_{2,2}t^2$.
    \item[(2.)] $x \mapsto x - \sqrt[3]{c_{2,2}}t^2, y \mapsto y$ yields the equation $y^2 = x^3 + t^3sx + a_{6,5}t^5s + t^4s^2$.
    \item[(3.)] We have $\Delta= -t^9s^3$ and by Tate's algorithm the fiber at $s=0$ is reducible if and only if $a_{6,5}=0$ in which case it has two components; so the RDP configuration on $X$ is $E_6 +A_1$ in this case.
    \end{enumerate}
    
    \item[---] $c_0 \neq 0$:
    \begin{enumerate}
    \item [(1.)] Rescaling $t$ and $s$, we can assume $c_0=1$. Then, we have $\Delta= -t^9 (c_2t-c_1^2t + c_1^3)$. Since $t \nmid c_1$, we can apply a substitution of the form $s \mapsto \mu s+ \lambda t $ for appropriate $\mu, \lambda \in k$ such that $s \mid \Delta$ and the coefficient of $-t^9s^3$ in $\Delta $ is $1$.
    \item [(2.)] $x \mapsto x + c_1t, y \mapsto y$ yields the equation $y^2=x^3 +t^2x^2+a_{6,5}t^5s + a_{6,4}t^4s^2 + t^3s^3$.
    \item[(3.)] We have $\Delta=-t^9s(a_{6,5}t^2 + a_{6,4}ts +s^2)$ and we see by Tate's algorithm that the RPD configuration on $X$ is $E_6+A_2$ if $a_{6,5}=a_{6,4}=0$, that it is $E_6 + A_1$ if $a_{6,4} \neq 0$ and $(a_{6,5}=0 \text{ or } a_{6,5}=a_{6,4}^2)$, and $E_6$ otherwise.
    \end{enumerate}
\end{itemize}

\subsubsection*{Lang's Type 6B (${\rm IV}^*$)}
$X$ is given by $y^2= x^3+ c_0t^2x^2  + d_0t^4x + c_2t^4$ with $t \nmid c_0, t \nmid c_2$ and $c_i, d_i \in k[t,s]$ homogeneous of degree $i$.
\begin{enumerate}
    \item[(1.)]Rescaling $t$ and $s$, we can assume $c_0=1$. Then, we have $\Delta= - t^{10}(c_2 -d_0^2 t^2 + d_0^3t^2)$. Since $t \nmid c_2$, we can apply a substitution of the form $s \mapsto \mu s+ \lambda t $ for appropriate $\mu, \lambda \in k$ such that $s \mid \Delta$ and the coefficient of $-t^{10}s^2$ in $\Delta $ is $1$.
    \item[(2.)] $x \mapsto x+ d_0t^2, y \mapsto y$ yields the equation $y^2= x^3 + t^2x^2 + a_{6,5}t^5s + t^4s^2$.
    \item[(3.)] We have $\Delta= -t^{10}s(a_{6,5}t +s)$ and we see by Tate's algorithm that the RDP configuration on $X$ is $E_6 + A_1$ if $a_{6,5}=0$, and $E_6$ otherwise.
\end{enumerate}

\subsubsection*{Lang's Type 6C (${\rm IV}^*$)}
$X$ is given by $y^2= x^3+  d_0t^4x + c_2t^4$ with $t \nmid d_0, t \nmid c_2$ and $c_i, d_i \in k[t,s]$ homogeneous of degree $i$.
\begin{enumerate}
    \item[(1.)] Using an automorphism of $\bbP^1$ we can assume that $d_0=1$ and $c_2= s^2 + c_{2,2}t^2$.
    \item[(2.)] $x \mapsto x+ \lambda t^2, y \mapsto y$ with $\lambda^3 + \lambda + c_{2,2}=0$ yields the equation $y^2= x^3 + t^4x + t^4$.
    \item[(3.)] Since $\Delta= -t^{12}$, $X$ has no other singularities apart from $E_6$.
\end{enumerate}

\subsubsection*{Lang's Type 7 (${\rm III}^*$)}
$X$ is given by $y^2= x^3+  c_0t^2x^2 + c_1t^3x+ d_1t^5$ with $t \nmid c_1$ and $c_i, d_i \in k[t,s]$ homogeneous of degree $i$.
From now on, let us distinguish the cases $c_0=0$ and $c_0 \neq 0$.
\begin{itemize}
    \item[---] $c_0 = 0$:
    \begin{enumerate}
    \item[(1.)] Since $t \nmid c_1$, we can apply an automorphism of $\bbP^1$ to assume $c_1=s$. Write $d_1= d_{1,0}s + d_{1,1}t$. Then, scaling $s \mapsto \lambda^3 s, t \mapsto \lambda^{-1}t$ for an appropriate $\lambda$, we can assume $d_{1,0}^3- d_{1,1} \in \{0,1\}$.
    \item[(2.)] $x \mapsto x - \sqrt[3]{d_{1,1}}t^2, y \mapsto y$ yields the equation $y^2 = x^3 + t^3sx + (d_{1,0}- \sqrt[3]{d_{1,1}})t^5s$.
    \item[(3.)] We have $\Delta= -t^9s^3$ and by Tate's algorithm the fiber at $s=0$ is reducible if and only if $d_{1,0}^3- d_{1,1}=0$ in which case it has two components; so the RDP configuration on $X$ is $E_7 +A_1$ in this case, and $E_7$ if $d_{1,0}^3- d_{1,1}=1$.
    \end{enumerate}

    \item[---] $c_0 \neq 0$:
    \begin{enumerate}
    \item [(1.)] Rescaling $t$ and $s$, we can assume $c_0=1$. Then, we have $\Delta= -t^9 (d_1 t^2 - c_1^2 t + c_1^3)$. Since $t \nmid c_1$, we can apply a substitution of the form $s \mapsto \mu s+ \lambda t $ for appropriate $\mu, \lambda \in k$ such that $s \mid \Delta$ and the coefficient of $-t^9s^3$ in $\Delta $ is $1$.
    \item [(2.)] $x \mapsto x + c_1t, y \mapsto y$ yields the equation $y^2=x^3 +t^2x^2+a_{6,5}t^5s - t^4s^2 + t^3s^3$.
    \item[(3.)] We have $\Delta=-t^9s(a_{6,5}t^2 -ts + s^2)$ and we see by Tate's algorithm that the RPD configuration on $X$ is $E_7+A_1$ if $a_{6,5} \in \{0,1\}$, and only $E_7$ otherwise.
    \end{enumerate}
\end{itemize}

\subsubsection*{Lang's Types 8A and 8B (${\rm II}^*$)}
These equations have been simplified by Lang in \cite{ExtremalCharpII} and they are as described in the second column of Table \ref{TableE6E7E8-char3}.

\subsubsection*{Quasi-elliptic surfaces: Ito's Types 3.3(1) (${\rm II}^*$) and 3.3(2) (${\rm IV}^*$)}
These equations have been simplified by Ito in \cite[Theorem 3.3]{QuasiEllipticChar3} and they are as described in the second column of Table \ref{TableE6E7E8-char3}.
\end{proof}

\begin{Example} \label{Ex dim def char 3}
Consider Lang's Type 7 with $c_0 \neq 0$ simplified as in the proof of Proposition \ref{Prop Equations-char3} and localized at the $E_7$ singularity at $t = x = y = 0$:
$$
f = y^2 - (x^3 + t^2 x^2 + a_{6,5}t^5 - t^4 + t^3) = 0.
$$
Applying Proposition \ref{Prop dimension(def)}, we have
\begin{eqnarray*}
T_f &=& k[t,x,y]_{(t,x,y)}/(f, tx^2 + a_{6,5}t^4 - t^3, xt^2, -y) \\
 &\cong& k[t,x]_{(t,x)}/(x^3 + a_{6,5}t^5 - t^4 + t^3, tx^2 + a_{6,5}t^4-t^3, xt^2) \\
 &\cong& k[t,x]_{(t,x)}/(x^3 + t^3, tx^2 - t^3, xt^2) =: R,
\end{eqnarray*}
where for the second isomorphism we have used that $0= t^2x^2 + a_{6,5}t^5 + t^4 = t^4(1 + a_{6,5}t)$, hence $t^4 = 0$, as $1+ a_{6,5}t$ is a unit.
Now, it is easy to check that $R$ is generated as a $k$-vector space by $1,x,x^2,x^3,t,tx,tx^2,tx^3,t^2$, hence $\dim T_f = 7$. Therefore, by Table \ref{TableNonTaut-char235}, the Artin coindex of this $E_7$ singularity is $1$.
\end{Example}

\subsection{Proof of Theorem \ref{Thm main} in Characteristic $2$} \label{Subsec Proof of Thm main char 2}
Assume $p = 2$. Let $\Gamma = \sum_i \Gamma_{i,n_i}^{k_i}$ be an RDP configuration containing a non-taut RDP and let $\Gamma' = \sum_{i} \Gamma_{i,n_i}$ be the associated configuration of $(-2)$-curves.
By Table \ref{TableNonTaut-char235}, we have $\Gamma_{i,n_i} \in \{D_n,E_6,E_7,E_8\}$ for some $i$.

First, observe that if $\Gamma \in \{E_7^0,D_6^0 + A_1, D_4^0 + 3A_1\}$, then $\Gamma$ occurs on a weak del Pezzo surface of degree $2$ by Proposition \ref{Prop RDP iff RDP deg 1}. Moreover, its associated $\Gamma'$ embeds into $E_8$ and $\Gamma' \not \in \{2A_3 + 2A_1, A_3 + 4A_1, 6A_1\}$, so Theorem \ref{Thm main} holds for these three exceptional cases.

Next, if $\Gamma \not \in  \{E_7^0,D_6^0 + A_1, D_4^0 + 3A_1\}$, then, by Proposition \ref{Prop RDP iff RDP deg 1}, $\Gamma$ occurs on an RDP del Pezzo surface if and only if it occurs on an RDP del Pezzo surface of degree $1$. Hence, by Lemma \ref{Lemma if gamma occurs, then conditions satisfied}, Theorem \ref{Thm main} holds for all $\Gamma$ such that $\Gamma'$ does not embed into $E_8$. Thus, we may assume that $\Gamma'$ embeds into $E_8$ and we note that $\Gamma' \not \in \{2A_3 + 2A_1, A_3 + 4A_1, 6A_1\}$, since $\Gamma$ contains a non-taut summand. 
On the other hand, it will follow from Proposition \ref{Prop Equations-char2} that every such $\Gamma'$ occurs on some weak del Pezzo surface of degree $1$. Thus, the following Proposition \ref{Prop Equations-char2} finishes the proof of Theorem \ref{Thm main}.

\begin{Proposition} \label{Prop Equations-char2} 
An RDP configuration $\Gamma$ containing a non-taut rational double point occurs on an RDP del Pezzo surface if and only if it occurs in Table \ref{TableD4D5-char2}, \ref{TableD6D7D8-char2}, or \ref{TableE6E7E8-char2}. 
Moreover, every RDP del Pezzo surface of degree $1$ containing a non-taut rational double point is given by an equation in one of these tables.
\end{Proposition}

\begin{proof}
This time, not all of the Weierstra{\ss} equations for $X \subseteq \bbP(1,1,2,3)$ in the classification of rational (quasi-)elliptic surfaces in \cite{EllipticChar2}, \cite{ExtremalCharpII} and \cite{QuasiEllipticChar2} are of the form (\hyperref[W2]{W2}) or (\hyperref[W2']{W2$'$}). Thus, to simplify them, we have to add a $0$th Step, before we can go on with our procedure as follows:
\begin{enumerate}
    \item[(0.)] Transform the Weierstra{\ss} equation into the form (\hyperref[W2]{W2}) or (\hyperref[W2']{W2$'$}).
    \item[(1.)] Carry out a substitution in $t$ and $s$ only.
    \item[(2.)] Apply an automorphism of $\bbP(1,1,2,3)$ over $\bbP^1$ preserving the form (\hyperref[W2]{W2}), or (\hyperref[W2']{W2$'$}), respectively, as in Subsection \ref{Subsub general Auto-char2}.
    \item[(3.)] Check for additional rational double points (\textit{e.g.} using Tate's algorithm (see Subsection \ref{Subsec Tate algo}) to determine the other reducible fibers of the underlying rational (quasi-)elliptic surface).
    \item[[(4.)] Determine the Artin coindices as described in Section \ref{Section Non-taut rational double points in positive characteristic}, \emph{e.g.} via Proposition \ref{Prop dimension(def)}.
 This will be left to the reader, but we will show how it works in Example \ref{Ex dim def char 2}.]
\end{enumerate}

\subsubsection*{Lang's Type 4A (${\rm I}_0^*$)}
$X$ is given by $y^2+ txy+ c_1t^2y= x^3+ d_1tx^2 + e_1t^3x + c_3t^3$ with $t \nmid c_1, t \nmid c_3$ and $c_i, d_i, e_i \in k[t,s]$ homogeneous of degree $i$. 
    \begin{enumerate}
    \item[(0.)] We want to transform the Weierstra{\ss} equation into one of the form (\hyperref[W2]{W2}). For this, send $x \mapsto c_1t$ to obtain the new equation
    $
        y^2 + txy= x^3 + d_1 tx^2 + e_2 t^2 x + c_3 t^3
    $
    with $t \nmid e_2$ and $t \nmid (e_2d_1 + c_3)$.
    \item[(1.)] For this new equation, we have $\Delta= t^8 (c_3t + e_2^2)$. Since $t \nmid e_2$, we can apply a substitution of the form $s \mapsto \mu s + \lambda t$ for appropriate $\mu, \lambda \in k$ such that $s \mid \Delta$ and the coefficient of $t^8s^4$ in $\Delta$ is $1$.
    \item[(2.)] We write $d_1= d_{1,0}s + d_{1,1}t$. Then, $x \mapsto x, y \mapsto y + \lambda t x + e_2 t$, where $\lambda$ is chosen in such a way that $\lambda^2 + \lambda = d_{1,1}$. This yields the equation 
    \begin{equation*}
        y^2 + txy= x^3 + a_{2,1}tsx^2 + a_{6,5}t^5s + a_{6,4}t^4s^2 + a_{6,3}t^3s^3 + t^2s^4 \text{ with } a_{2,1}+ a_{6,3} \neq 0 \text{.}
    \end{equation*}
    \item[(3.)] We have $\Delta= t^8s (a_{6,5}t^3 + a_{6,4}t^2s + a_{6,3}t s^2 + s^3)$. By Tate's algorithm we see that the RDP configuration on $X$ is $D_4 + A_3$ if $a_{6,5}= a_{6,4}= a_{6,3}=0$, that it is $D_4 + A_2$ if $(a_{6,5}= a_{6,4}=0 \text{ and } a_{6,3} \neq 0)$ or $(a_{6,3}^2 = a_{6,4} \text{ and } a_{6,3}^3= a_{6,5} \neq 0)$, that it is $D_4 + 2A_1$ if $a_{6,5}=a_{6,3}=0$ and $a_{6,4} \neq 0$, that it is $D_4 + A_1$ if $(a_{6,5}=0 \text{ and } a_{6,3} \neq 0)$ or $a_{6,5}= a_{6,4}a_{6,3} \neq 0$, and that it is only $D_4$ otherwise.
    \end{enumerate}

\subsubsection*{Lang's Type 4B (${\rm I}_0^*$)}
$X$ is given by $y^2+ txy+ c_0 t^3y= x^3+ d_1tx^2 + e_1t^3x + c_3t^3$ with $t \nmid c_3$ and $c_i, d_i, e_i \in k[t,s]$ homogeneous of degree $i$. 
    \begin{enumerate}
    \item[(0.)] We want to transform the Weierstra{\ss} equation into one of the form (\hyperref[W2]{W2}). For this, send $x \mapsto c_0t^2$ to obtain the new equation
    $
        y^2 + txy= x^3 + d_1 tx^2 + e_1 t^3x + c_3 t^3
    $
    with $t \nmid c_3$.
    \item[(1.)] For this new equation, we have $\Delta= t^9 (c_3+ e_1^2 t)$. Since $t \nmid c_3$, we can apply a substitution of the form $s \mapsto \mu s + \lambda t$ for appropriate $\mu, \lambda \in k$ such that $s \mid \Delta$ and the coefficient of $t^9s^3$ in $\Delta$ is $1$.
    \item[(2.)] We write $d_1= d_{1,0}s + d_{1,1}t$. Then, $x \mapsto x, y \mapsto y + \lambda t x + e_1 t^2$, where $\lambda$ is chosen in such a way that $\lambda^2 + \lambda = d_{1,1}$. This yields the equation 
    $
        y^2 + txy= x^3 + a_{2,1}tsx^2 + a_{6,5}t^5s + a_{6,4}t^4s^2 + t^3s^3 \text{.}
    $
    \item[(3.)] We have $\Delta= t^9s (a_{6,5}t^2 + a_{6,4}ts + s^2)$. By Tate's algorithm we see that the RDP configuration on $X$ is $D_4 + A_2$ if $a_{6,5}= a_{6,4}=0$, that it is $D_4 + A_1$ if $(a_{6,5}=0 \text{ and } a_{6,4} \neq 0)$ or $(a_{6,5} \neq 0 \text{ and } a_{6,4}= 0)$, and that it is only $D_4$ otherwise.
    \end{enumerate}

\subsubsection*{Lang's Type 5A (${\rm I}_1^*$)}
$X$ is given by $y^2+ txy+ c_1t^2y= x^3+ d_1tx^2 + e_1t^3x + c_2t^4$ with $t \nmid c_1, t \nmid d_1$ and $c_i, d_i, e_i \in k[t,s]$ homogeneous of degree $i$. 
    \begin{enumerate}
    \item[(0.)] We want to transform the Weierstra{\ss} equation into one of the form (\hyperref[W2]{W2}). For this, send $x \mapsto c_1t$ to obtain the new equation
    $
        y^2 + txy= x^3 + d_1 tx^2 + e_2 t^2 x + c_3 t^3
    $
    with $t \nmid e_2$ and $t \mid  (e_2d_1 + c_3)$.
    \item[(1.)] For this new equation, we have $\Delta= t^8 (c_3t + e_2^2)$. Since $t \nmid e_2$, we can apply a substitution of the form $s \mapsto \mu s + \lambda t$ for appropriate $\mu, \lambda \in k$ such that $s \mid \Delta$ and the coefficient of $t^8s^4$ in $\Delta$ is $1$.
    \item[(2.)] We write $d_1= d_{1,0}s + d_{1,1}t$. Then, $x \mapsto x, y \mapsto y + \lambda t x + e_2 t$, where $\lambda$ is chosen in such a way that $\lambda^2 + \lambda = d_{1,1}$. This yields the equation 
    \begin{equation*}
        y^2 + txy= x^3 + a_{2,1}tsx^2 + a_{6,5}t^5s + a_{6,4}t^4s^2 + a_{6,3}t^3s^3 + t^2s^4 \text{ with } a_{2,1} + a_{6,3}=0 \text{.}
    \end{equation*}
    \item[(3.)] We have $\Delta= t^8s (a_{6,5}t^3 + a_{6,4}t^2s + a_{6,3}t s^2 + s^3)$. By Tate's algorithm we see that the RDP configuration on $X$ is $D_5 + A_3$ if $a_{6,5}= a_{6,4}= a_{6,3}=0$, that it is $D_5 + A_2$ if $(a_{6,5}= a_{6,4}=0 \text{ and } a_{6,3} \neq 0)$ or $(a_{6,3}^2 = a_{6,4} \text{ and } a_{6,3}^3= a_{6,5} \neq 0)$, that it is $D_5 + 2A_1$ if $a_{6,5}=a_{6,3}=0$ and $a_{6,4} \neq 0$, that it is $D_5 + A_1$ if $(a_{6,5}=0 \text{ and } a_{6,3} \neq 0)$ or $a_{6,5}= a_{6,4}a_{6,3} \neq 0$, and that it is only $D_5$ otherwise.
    \end{enumerate}

    \subsubsection*{Lang's Type 5B (${\rm I}_2^*$)}
    First, we note that there seems to be a typo in \cite[Case 5B., p.5825]{EllipticChar2} in the sense that the fiber type should be (${\rm I}_2^*$) instead of (${\rm I}_1^*$).
$X$ is given by $y^2+ txy+ c_0 t^3y= x^3+ d_1tx^2 + e_1t^3x + c_1t^5$ with $t \nmid d_1, t \nmid e_1$ and $c_i, d_i, e_i \in k[t,s]$ homogeneous of degree $i$. 
    \begin{enumerate}
    \item[(0.)] We want to transform the Weierstra{\ss} equation into one of the form (\hyperref[W2]{W2}). For this, send $x \mapsto c_0t^2$ to obtain the new equation
    $
        y^2 + txy= x^3 + d_1 tx^2 + e_1 t^3x + c_1 t^5
    $
    with $t \nmid d_1$ and $t \nmid e_1$.
    \item[(1.)] For this new equation, we have $\Delta= t^{10} (e_1^2 +c_1t)$. Since $t \nmid e_1$, we can apply a substitution of the form $s \mapsto \mu s + \lambda t$ for appropriate $\mu, \lambda \in k$ such that $s \mid \Delta$ and the coefficient of $t^{10}s^2$ in $\Delta$ is $1$.
    \item[(2.)] We write $d_1= d_{1,0}s + d_{1,1}t$, where $d_{1,0} \neq 0$. Then, $x \mapsto x, y \mapsto y + \lambda t x + e_1 t^2$, where $\lambda$ is chosen in such a way that $\lambda^2 + \lambda = d_{1,1}$. This yields the equation 
    $
        y^2 + txy= x^3 + a_{2,1}tsx^2 + a_{6,5}t^5s + t^4s^2 
    $ with $a_{2,1} \neq 0$.
    \item[(3.)] We have $\Delta= t^{10}s (a_{6,5}t + s)$. By Tate's algorithm we see that the RDP configuration on $X$ is $D_6 + A_1$ if $a_{6,5}=0$, and $D_6$ otherwise.
    \end{enumerate}

        \subsubsection*{Lang's Type 5C (${\rm I}_3^*$)}
$X$ is given by $y^2+ txy+ c_0 t^3y= x^3+ d_1tx^2 + e_0t^4x + d_0t^6$ with $t \nmid c_0, t \nmid d_1$ and $c_i, d_i, e_i \in k[t,s]$ homogeneous of degree $i$. 
    \begin{enumerate}
    \item[(0.)] We want to transform the Weierstra{\ss} equation into one of the form (\hyperref[W2]{W2}). For this, send $x \mapsto c_0t^2$ to obtain the new equation
    $
        y^2 + txy= x^3 + d_1 tx^2 + e_0t^4x + c_1 t^5
    $
    with $t \nmid d_1$ and $t \nmid c_1$.
    \item[(1.)] For this new equation, we have $\Delta= t^{11} (c_1+ e_0^2t)$. Since $t \nmid c_1$, we can apply a substitution of the form $s \mapsto \mu s + \lambda t$ for appropriate $\mu, \lambda \in k$ such that $s \mid \Delta$, \emph{i.e.}, $(c_1+ e_0^2t)=s$.
    \item[(2.)] We write $d_1= d_{1,0}s + d_{1,1}t$, where $d_{1,0} \neq 0$. Then, $x \mapsto x, y \mapsto y + \lambda t x + e_0t^3$, where $\lambda$ is chosen in such a way that $\lambda^2 + \lambda = d_{1,1}$. This yields the equation 
    $
        y^2 + txy= x^3 + a_{2,1}tsx^2 + t^5s 
    $ with $a_{2,1} \neq 0$.
    \item[(3.)] We have $\Delta= t^{11}s$ and the RDP configuration on $X$ is $D_7$.
    \end{enumerate}
    
    \subsubsection*{Lang's Type 5D (${\rm I}_4^*$)}
This equation has been simplified by Lang in \cite{ExtremalCharpII} and it is as described in the second column of Table \ref{TableD6D7D8-char2}.

\subsubsection*{Lang's Type 6 (${\rm IV}^*$)}
$X$ is given by $y^2+ txy+ c_1t^2y= x^3+ d_0t^2x^2 + e_1t^3x + c_2t^4$ with $t \nmid c_1$ and $c_i, d_i, e_i \in k[t,s]$ homogeneous of degree $i$. 
    \begin{enumerate}
    \item[(0.)] We want to transform the Weierstra{\ss} equation into one of the form (\hyperref[W2]{W2}). For this, send $x \mapsto c_1t$ to obtain the new equation
    $
        y^2 + txy= x^3 + d_1 tx^2 + e_2t^2x + c_3 t^3
    $
    with $t \nmid e_2, t \mid (d_1^2 + e_2)$ and $t \mid (d_1e_2 + c_3)$.
    \item[(1.)] For this new equation, we have $\Delta= t^{8} (c_3t + e_2^2)$. Since $t \nmid e_2$, we can apply a substitution of the form $s \mapsto \mu s + \lambda t$ for appropriate $\mu, \lambda \in k$ such that $s \mid \Delta$ and the coefficient of $t^{8}s^4$ in $\Delta$ is $1$.
    \item[(2.)] We write $d_1= d_{1,0}s + d_{1,1}t$, where $d_{1,0} \neq 0$. Then, $x \mapsto x, y \mapsto y + \lambda t x + e_2t$, where $\lambda$ is chosen in such a way that $\lambda^2 + \lambda = d_{1,1}$. This yields the equation 
    \begin{equation*}
        y^2 + txy= x^3 + tsx^2 + a_{6,5}t^5s + a_{6,4}t^4s^2 + t^3s^3 + t^2s^4 \text{.}
    \end{equation*}
    \item[(3.)] We have $\Delta= t^{8}s(a_{6,5}t^3 + a_{6,4}t^2s + ts^2 + s^3)$. By Tate's algorithm we see that the RDP configuration on $X$ is $E_6 + A_2$ if $a_{6,5}=a_{6,4} \in \{0,1\}$, that it is $E_6 + A_1$ if ($a_{6,5}=0$ and $a_{6,4} \neq 0$) or $a_{6,5}= a_{6,4} \not\in \{0,1\}$, and that it is only $E_6$ otherwise.
    \end{enumerate}

    \subsubsection*{Lang's Type 7 (${\rm III}^*$)}
$X$ is given by $y^2+ txy+ c_0t^3y= x^3+ d_0t^2x^2 + e_1t^3x + f_1t^5$ with $t \nmid e_1$ and $c_i, d_i, e_i, f_i \in k[t,s]$ homogeneous of degree $i$. 
    \begin{enumerate}
    \item[(0.)] We want to transform the Weierstra{\ss} equation into one of the form (\hyperref[W2]{W2}). For this, send $x \mapsto c_0t^2$ to obtain the new equation
    $
        y^2 + txy= x^3 + d_0t^2x^2 + e_1t^3x + f_1t^5
    $
    with $t \nmid e_1$.
    \item[(1.)] For this new equation, we have $\Delta= t^{10} (e_1^2 + f_1t)$. Since $t \nmid e_1$, we can apply a substitution of the form $s \mapsto \mu s + \lambda t$ for appropriate $\mu, \lambda \in k$ such that $s \mid \Delta$ and the coefficient of $t^{10}s^2$ in $\Delta$ is $1$.
    \item[(2.)] Then, $x \mapsto x, y \mapsto y + \lambda t x + e_1 t^2$, where $\lambda$ is chosen in such a way that $\lambda^2 + \lambda = d_{0}$. This yields the equation 
    $
        y^2 + txy= x^3 +  a_{6,5}t^5s +t^4s^2 
    $.
    \item[(3.)] We have $\Delta= t^{10}s(a_{6,5}t + s)$. By Tate's algorithm we see that the RDP configuration on $X$ is $E_7 + A_1$ if $a_{6,5}=0$, and $E_7$ otherwise.
    \end{enumerate}
    
    \subsubsection*{Lang's Type 8 (${\rm II}^*$)}
This equation has been simplified by Lang in \cite{ExtremalCharpII} and it is as described in the second column of Table \ref{TableD6D7D8-char2}.

\subsubsection*{Lang's Type 12A (${\rm I}_0^*$)}
$X$ is given by $y^2+ c_1 t^2y= x^3+ d_1tx^2 + e_1 t^3x + c_3t^3$ with $t \nmid c_1, t \nmid c_3$ and $c_i, d_i, e_i \in k[t,s]$ homogeneous of degree $i$. 
    \begin{enumerate}
    \item[(0.)] The Weierstra{\ss} equation is already of the form (\hyperref[W2']{W2$'$}). 
    \item[(1.)] We have $\Delta= t^{8} c_1^4$. Since $t \nmid c_1$, we can apply a substitution of the form $s \mapsto \mu s + \lambda t$ for appropriate $\mu, \lambda \in k$ such that $c_1=s$. 
    We write $c_3= c_{3,0}s^3 + c_{3,1}ts^2 + c_{3,2}t^2s + c_{3,3}t^3$.
    Then, since $t \nmid c_3$, scaling $s \to \lambda^2 s, t \mapsto \lambda^{-1}t$ for an appropriate $\lambda$ yields $c_{3,0}=1$.
    \item[(2.)] Let us write $d_1= d_{1,0}s + d_{1,1}t$ and $e_1= e_{1,0}s + e_{1,1}t$. Then, we carry out the substitution 
    \begin{equation*}
    x \mapsto x + \sqrt{e_{1,1}}t^2, \hspace{3mm} y \mapsto y + \sqrt{d_{1,1} + \sqrt{e_{1,1}}}tx + \lambda t^2s + (c_{3,2} + e_{1,0} \sqrt{e_{1,1}}+ d_{1,0}e_{1,1})t^3 \text{,}
    \end{equation*}
    where $\lambda^2 + \lambda = c_{3,1}$. This yields
    $
        y^2 + t^2s y= x^3+ a_{2,1}ts x^2 + a_{4,3}t^3s x + a_{6,6}t^6 + t^3s^3 
    $.
    \item[(3.)] We have $\Delta= t^8s^4$. By Tate's algorithm we see that the RDP configuration on $X$ is $D_4 + A_2$ if $a_{6,6}=a_{4,3} = 0$, that it is $D_4 + A_1$ if $a_{6,6}=0 \text{ and } a_{4,3} \neq 0$, and that it is $D_4$ otherwise.
    \end{enumerate}

    \subsubsection*{Lang's Type 12B (${\rm I}_0^*$)}
$X$ is given by $y^2+ c_0 t^3y= x^3+ d_1tx^2 + c_1 t^3x + c_3t^3$ with $t \nmid c_0, t \nmid c_3$ and $c_i, d_i \in k[t,s]$ homogeneous of degree $i$. For the simplification of this equation, we will not follow the procedure described in the beginning of the proof, but perform the substitutions in a different order.
  \begin{enumerate}
    \item[(a.)] First, applying $x \mapsto x + d_1t$ and then, scaling $x \mapsto \lambda^3 x, y \mapsto \lambda^2 y $ for an appropriate $\lambda$ yields the new equation $y^2+ t^3y= x^3 + c_2 t^2x + c_3t^3$ with $t \nmid c_3$.
    \item[(b.)] Since $t \nmid c_3$, we can apply a substitution of the form $s \mapsto \mu s + \lambda t$ for appropriate $\mu, \lambda \in k$ to obtain a $c_3$ of the form $c_3= s^3 + a^2s^2t + a s t^2+ c_{3,3}t^3$ for some $a \in k$.
    \item[(c.)] Finally, $y \mapsto y + ast^2 + \mu t^3$ with $\mu^2 + \mu= c_{3,3}$ yields 
    \begin{equation*}
    y^2 + t^3y= x^3 + (a_{4,2}s^2 + a_{4,3}ts + a_{4,4}t^2)t^2x + s^3t^3 \text{.}
    \end{equation*}
    \item[(3.)] Since $\Delta= t^{12}$, there are no other reducible fibers and the RDP configuration on $X$ is $D_4$.
    \end{enumerate}

\subsubsection*{Lang's Type 13A (${\rm I}_1^*$)}
$X$ is given by $y^2+ c_1 t^2y= x^3+ d_1tx^2 + e_1 t^3x + d_2t^4$ with $t \nmid c_1, t \nmid d_1$ and $c_i, d_i, e_i \in k[t,s]$ homogeneous of degree $i$. 
    \begin{enumerate}
    \item[(1.)] We have $\Delta= t^{8} c_1^4$. Since $t \nmid c_1$, we can apply a substitution of the form $s \mapsto \mu s + \lambda t$ for appropriate $\mu, \lambda \in k$ such that $c_1=s$. Then, since $t \nmid d_1$, we can scale $s$ and $t$ such that $d_1= s + d_{1,1}t$.
    \item[(2.)] Let us write $e_1= e_{1,0}s + e_{1,1}t$ and $d_2= d_{2,0}s^2 + d_{2,1}ts + d_{2,2}t^2$. The substitution 
    $
    x \mapsto x + \sqrt{e_{1,1}}t^2$,  
    $y \mapsto y + e_{1,0}tx 
    + \lambda t^2s 
    + \sqrt{d_{1,1}e_{1,1} + d_{2,2}} t^3 
    $
    with $\lambda^2 + \lambda = d_{2,0}$ yields the new equation
    \mbox{$
       y^2 + t^2s y= x^3 + (a_{2,2}t^2 + ts) x^2 + a_{6,5}t^5s
    $.}
    \item[(3.)] We have $\Delta= t^8s^4$. By Tate's algorithm we see that the RDP configuration on $X$ is $D_5 + A_2$ if $a_{6,5}=a_{2,2} = 0$, that it is $D_5 + A_1$ if $a_{6,5}=0 \text{ and } a_{2,2} \neq 0$, and that it is $D_5$ otherwise.
    \end{enumerate}

    \subsubsection*{Lang's Type 13B (${\rm I}_2^*$)}
$X$ is given by $y^2+ c_0 t^3y= x^3+ d_1tx^2 + e_1 t^3x + f_1 t^5$ with $t \nmid c_0, t \nmid e_1$ and $c_i, d_i, e_i, f_i \in k[t,s]$ homogeneous of degree $i$. 
    \begin{enumerate}
    \item[(1.)] Since $t \nmid e_1$, we can apply a substitution of the form $s \mapsto \mu s + \lambda t$ for appropriate $\mu, \lambda \in k$ such that $e_1=s$. Then, since $t \nmid c_0$, we can scale $s$ and $t$ such that $c_0= 1$.
    \item[(2.)] Let $d_1= d_{1,0}s + d_{1,1}t$ and $f_1= f_{1,0}s + f_{1,1}t$. The substitution 
    $
    x \mapsto x + \lambda t^2,  
    y \mapsto y + \lambda^2 tx 
    + \mu t^3 $,
    where $\lambda$ and $\mu$ are chosen such that $ d_{1,0} \lambda^2 + \lambda= f_{1,0}$ and $\mu^2 + \mu= \lambda^3 + d_{1,1} \lambda^2 + f_{1,1}$, yields the new equation
    \mbox{$
       y^2 + t^3y= x^3 + (a_{2,2}t^2 + a_{2,1}ts) x^2 + t^3sx
    $.}
    \item[(3.)] Since $\Delta= t^{12}$, the RDP configuration on $X$ is $D_6$.
    \end{enumerate}

\subsubsection*{Lang's Type 13C (${\rm I}_3^*$)}
$X$ is given by $y^2+ c_0 t^3y= x^3+ d_1tx^2 + d_0 t^4 x + e_0t^6$ with $t \nmid c_0, t \nmid d_1$ and $c_i, d_i, e_i \in k[t,s]$ homogeneous of degree $i$. 
    \begin{enumerate}
    \item[(1.)] Since $t \nmid d_1$, we can apply a substitution of the form $s \mapsto \mu s + \lambda t$ for appropriate $\mu, \lambda \in k$ such that $d_1=s$. Then, since $t \nmid c_0$, we can scale $s$ and $t$ such that $c_0= 1$. Further, we apply the $\bbP^1$-automorphism $s \mapsto s + d_0^2 t$ to obtain $y^2+  t^3y= x^3+ tsx^2 + d_0^2 t^2 x^2 + d_0 t^4 x + e_0t^6$.
    \item[(2.)] Then the substitution
    $
    y \mapsto y + d_0 tx + \lambda t^3 $
    with $\lambda^2 + \lambda= e_0$ yields $y^2+  t^3y= x^3+ tsx^2  $.
    \item[(3.)] Since $\Delta= t^{12}$, the RDP configuration on $X$ is $D_7$.
    \end{enumerate}

    \subsubsection*{Lang's Type 14 (${\rm IV}^*$)}
$X$ is given by $y^2+ c_1 t^2y= x^3+ d_0t^2x^2 + e_1 t^3x + d_2t^4$ with $t \nmid c_1$ and $c_i, d_i, e_i \in k[t,s]$ homogeneous of degree $i$. 
    \begin{enumerate}
    \item[(1.)] Since $t \nmid c_1$, we can apply a substitution of the form $s \mapsto \mu s + \lambda t$ for appropriate $\mu, \lambda \in k$ such that $c_1=s$. 
    \item[(2.)] Let us write $e_1= e_{1,0}s + e_{1,1}t$ and $d_2= d_{2,0}s^2 + d_{2,1}ts + d_{2,2}t^2$. The substitution 
    \mbox{$
    x \mapsto x + \sqrt{e_{1,1}} t^2$,  }
    $y \mapsto y + e_{1,0}tx 
    + \lambda t^2s 
    + \sqrt{d_0 e_{1,1} + d_{2,2}} t^3 
    $
    with $\lambda^2 + \lambda = d_{2,0}$ yields the new equation
    $
       y^2 + t^2s y$ $= x^3 + a_{2,2}t^2 x^2 + a_{6,5}t^5s \text{.}
    $
    \item[(3.)] We have $\Delta= t^8s^4$. By Tate's algorithm we see that the RDP configuration on $X$ is $E_6 + A_2$ if $a_{6,5}=  a_{2,2} = 0$, that it is $E_6 + A_1$ if $a_{6,5}= 0 \text{ and } a_{2,2} \neq 0$, and that it is $E_6$ otherwise.
    \end{enumerate}

    \subsubsection*{Lang's Type 15 (${\rm III}^*$)}
$X$ is given by $y^2+ c_0 t^3y= x^3+ d_0 t^2 x^2 + e_1 t^3x + d_1 t^5$ with $t \nmid c_0, t \nmid e_1$ and $c_i, d_i, e_i \in k[t,s]$ homogeneous of degree $i$. 
    \begin{enumerate}
    \item[(1.)] Since $t \nmid e_1$, we can apply a substitution of the form $s \mapsto \mu s + \lambda t$ for appropriate $\mu, \lambda \in k$ such that $e_1=s$. Then, since $t \nmid c_0$, we can scale $s$ and $t$ such that $c_0= 1$. 
    Let us write $d_1= d_{1,0}s + d_{1,1}t$ and apply the $\bbP^1$-automorphism $s \mapsto s + (d_{1,0}^2 + \sqrt{d_0 + d_{1,0}}) t$ to obtain $y^2+ t^3y= x^3+ d_0 t^2 x^2 + t^3sx + (d_{1,0}^2 + \sqrt{d_0 + d_{1,0}})t^4x + d_{1,0}t^5s + d_{1,0}(d_{1,0}^2 + \sqrt{d_0 + d_{1,0}})t^6 + d_{1,1}t^6 $.
    \item[(2.)] Finally, the substitution
    $ x \mapsto x + d_{1,0}t^2, 
    y \mapsto y + \sqrt{d_{0} + d_{1,0}} tx + \lambda t^3 $
    with $\lambda^2 + \lambda= d_{1,0}^3 + d_0d_{1,0}^2 + d_{1,1}$ yields the simplified equation 
    $y^2+  t^3y= x^3+ t^3sx  $.
    \item[(3.)] Since $\Delta= t^{12}$, the RDP configuration on $X$ is $E_7$.
    \end{enumerate}

    \subsubsection*{Lang's Type 16 (${\rm II}^*$)}
This equation has been simplified by Lang in \cite{ExtremalCharpII} and it is as described in the second column of Table \ref{TableD6D7D8-char2}.

\subsubsection*{Ito's Types 5.2(a) (${\rm II}^*$), 5.2(b) (${\rm I}_4^*$), 5.2(c) (${\rm III}^*$), 5.2(d) (${\rm 2I}_0^*$), 5.2(e) (${\rm I}_2^*$) and 5.2(f) (${\rm I}_0^*$)}
    Note that $X$ is a quasi-elliptic surface.
These equations have been simplified by Ito in \cite[Theorem 5.2.]{QuasiEllipticChar2} and they are as described in the second column of Table \ref{TableE6E7E8-char3}, where we only simplified the equation for 5.2.(e) in order to put the $D_6$ singularity at $(t,x,y)=(0,0,0)$.
\end{proof}

\begin{Example} \label{Ex dim def char 2}

Consider Lang's Type 5B. simplified as in the proof of Proposition \ref{Prop Equations-char2} and localized at the $D_6$ singularity at $t = x = y = 0$:
$$
f = y^2 +txy + (x^3 + a_{2,1}t x^2 + a_{6,5}t^5 +t^4) = 0.
$$
Applying Proposition \ref{Prop dimension(def)}, we have
\begin{eqnarray*}
T_f &=& k[t,x,y]_{(t,x,y)}/(f, a_{2,1}x^2 + a_{6,5}t^4+xy, x^2 + ty, tx) \\
 &\cong& k[t,x,y]_{(t,x,y)}/(y^2 + a_{6,5}t^5 + t^4 , a_{2,1}x^2 + a_{6,5}t^4+xy, x^2 + ty, tx) \\
 &\cong& k[t,x,y]_{(t,x,y)}/(y^2 +  t^4 , a_{2,1}x^2 + a_{6,5}t^4+xy, x^2 + ty, tx) =: R,
\end{eqnarray*}
where for the first isomorphism we have used that $x^3= xty=0$ and for the second isomorphism we have used that $a_{6,5}t^5= a_{2,1}x^2t + xyt=0$.
Now, it is easy to check that $R$ is generated as a $k$-vector space by $1,x,y,t,x^2,y^2,t^2,t^3$, hence $\dim T_f = 8$. Therefore, by Table \ref{TableNonTaut-char235}, the Artin coindex of this $D_6$ singularity is $2$.

\end{Example}

\begin{Remark}\label{Rem Artin coindices which do not occur}
Using our results, it is straightforward to list all RDP configurations $\Gamma$ such that the associated configuration of $(-2)$-curves $\Gamma'$ occurs on a weak del Pezzo surface, but $\Gamma$ itself does not occur on any RDP del Pezzo surface. By Theorem \ref{Thm main}, this phenomenon happens only in characteristic $2$ and there precisely if
\begin{eqnarray*}
\Gamma &\in& \{E_8^1,E_8^2, E_7^1 + A_1, E_7^1, E_7^2 + A_1, D_8^1, D_8^2, D_7^0, D_6^0, D_6^1 + 2A_1, D_6^1 + A_1, D_6^2 + 2A_1, \\ && D_5^0 + A_3, D_5^0 + 2A_1, D_4^0+D_4^1, D_4^0 + A_3, D_4^1+D_4^1, D_4^1 + 4A_1, D_4^1 + 3A_1\}\text{.}
\end{eqnarray*}
It would be interesting to find an abstract reason for the non-existence of those Artin coindices on RDP del Pezzo surfaces.
\end{Remark}

\appendix
\section*{Appendix}\addcontentsline{toc}{section}{Appendix}

In the following tables, we list all possible RDP configurations containing a non-taut rational double point (it occurs in characteristic 2, 3, and 5).

\pagestyle{empty}

\newpage

\newgeometry{top=20 mm, bottom=8mm, left=0mm, right=0mm}

\begin{landscape}
\begin{table}[h] \renewcommand{\arraystretch}{1.35}
 $$
 \begin{array}{|c c c|c|c|c|c|c|} 
 \hline
\multicolumn{3}{|c|}{
\begin{array}{c}
    \text{RDP}\\
    \text{configuration}
\end{array}
}
& \begin{array}{c}
    \text{\hspace{7mm} Weierstra{\ss} equation of } X \text{ in } \bbP(1,1,2,3) \text{\hspace{7mm}} \\ 
    \hdashline
     \text{condition for extra RDPs}
\end{array}
& \Delta
& j
& \scriptstyle{\begin{array}{c}
\text{Lang's /}\\
\text{Ito's type}
\end{array} }
& \begin{array}{c}
\text{ell /}\\
\text{q-ell}
\end{array}
\\

\hline  \hline
%& 
\multicolumn{8}{|c|}{\bf{D_4}}
\\
\hline

%12B
\multicolumn{3}{|c|}{D_4^0}
& y^2 + t^3y= x^3 + (a_{4,2}s^2 + a_{4,3}ts + a_{4,4}t^2)t^2x + s^3t^3
& t^{12}
& 0
& \text{12B} \label{Tab12B}
& \text{ell}
\\
\hline

%12A
D_4^0
&
&
& y^2 + t^2s y= x^3+ a_{2,1}ts x^2 + a_{4,3}t^3s x + a_{6,6}t^6 + t^3s^3
& t^8s^4
& 0
& \text{12A} \label{Tab12A}
& \text{ell}
\\
\hdashline

& +
& A_1
& \text{if } a_{6,6}=0 \text{ and } a_{4,3} \neq 0
& 
& 
& \text{12A 10A} \label{Tab12A10A}
& \text{ell}
\\

& +
& A_2
& \text{if } a_{6,6}=a_{4,3} = 0
& 
& 
& \text{12A 11} \label{Tab12A11}
& \text{ell}
\\
\hline

%exception C
D_{4}^0
& + 
& 3A_1
& \multicolumn{5}{c|}{\text{occurs only in degree }$2$ \text{ (see Proposition \ref{Prop RDP iff RDP deg 1} (\hyperref[Exception C]{C.}))}}
\\
\hline

%5.2.(f)
D_4^0
& +
& 4A_1
& 
y^2 = x^3 + (t^3s+ a_{4,2}t^2s^2 +ts^3 )x \text{ with } a_{4,2} \neq 0
& 0
& 
&\text{5.2.(f)} \label{Tab52(f)}
& \text{q-ell}
\\
\hline

%5.2.(d)
D_4^0
& +
& D_4^0
& 
y^2 = x^3 + a_{4,2}t^2s^2x + t^3s^3
& 0
& 
&\text{5.2.(d)} \label{Tab52(d)}
& \text{q-ell}
\\
\hline

\hline

%4B.
D_4^1
&
&
& y^2 + txy= x^3 + a_{2,1}tsx^2 + a_{6,5}t^5s + a_{6,4}t^4s^2 + t^3s^3
& t^9s (a_{6,5}t^2 + a_{6,4}ts + s^2)
& \frac{t^{12}}{\Delta}
& \text{4B.} \label{Tab4B-D41}
& \text{ell}
\\
\hdashline

& +
& A_1
& \text{if } (a_{6,5}=0 \text{ and } a_{6,4} \neq 0) \text{ or } (a_{6,5} \neq 0 \text{ and } a_{6,4}= 0)
& 
& 
& \text{4B. 2.} \label{Tab4B2}
& \text{ell}
\\

& +
& A_2
& \text{if } a_{6,5}= a_{6,4}=0
& 
& 
& \text{4B. 3.} \label{Tab4B3}
& \text{ell}
\\
\hline

%4A.
D_4^1
&
&
&
\begin{array}{c}
     y^2 + txy= x^3 + a_{2,1}tsx^2 + a_{6,5}t^5s + a_{6,4}t^4s^2 + a_{6,3}t^3s^3 + t^2s^4  \\
     \text{with } a_{2,1}+ a_{6,3} \neq 0 
\end{array}
& t^8s (a_{6,5}t^3 + a_{6,4}t^2s + a_{6,3}t s^2 + s^3)
& \frac{t^{12}}{\Delta}
& \text{4A.} \label{Tab4A-D41}
& \text{ell}
\\
\hdashline

& +
& A_1
& \text{if } (a_{6,5}=0 \text{ and } a_{6,3} \neq 0) \text{ or } a_{6,5}= a_{6,4}a_{6,3} \neq 0
& 
& 
& \text{4A. 2.} \label{Tab4A2}
& \text{ell}
\\

& +
& 2A_1
& \text{if } a_{6,5}=a_{6,3}=0 \text{ and } a_{6,4} \neq 0
& 
& 
& \text{4A. 4.} \label{Tab4A4}
& \text{ell}
\\

& +
& A_2
& \text{if } (a_{6,5}= a_{6,4}=0 \text{ and } a_{6,3} \neq 0) \text{ or } (a_{6,3}^2 = a_{6,4} \text{ and } a_{6,3}^3= a_{6,5} \neq 0)
& 
& 
& \text{4A. 3.} \label{Tab4A3}
& \text{ell}
\\

& +
& A_3
& \text{if }  a_{6,5}= a_{6,4}= a_{6,3}=0
& 
& 
& \text{4A. 5.} \label{Tab4A5}
& \text{ell}
\\

\hline  \hline
%& 
\multicolumn{8}{|c|}{\bf{D_5}}
\\
\hline

%13A
D_5^0
&
&
& y^2 + t^2s y= x^3 + (a_{2,2}t^2 + ts) x^2 + a_{6,5}t^5s
& t^{8}s^4
& 0
& \text{13A} \label{Tab13A-D50}
& \text{ell}
\\
\hdashline

& +
& A_1
& \text{if }   a_{6,5}=0 \text{ and } a_{2,2} \neq 0
& 
& 
& \text{13A 10A} \label{Tab13A10A}
& \text{ell}
\\

& +
& A_2
& \text{if }   a_{6,5}=a_{2,2} = 0
& 
& 
& \text{13A 11} \label{Tab13A11}
& \text{ell}
\\
\hline

\hline

%5A.
D_5^1
&
&
& \begin{array}{c}
y^2 + txy= x^3 + a_{2,1}tsx^2 + a_{6,5}t^5s + a_{6,4}t^4s^2 + a_{6,3}t^3s^3 + t^2s^4 \\
\text{with } a_{2,1} = a_{6,3}
\end{array}
& 
t^8s (a_{6,5}t^3 + a_{6,4}t^2s + a_{6,3}t s^2 + s^3)
& \frac{t^{12}}{\Delta}
& \text{5A.} \label{Tab5A-D51}
& \text{ell}
\\
\hdashline

& +
& A_1
& \text{if } (a_{6,5}=0 \text{ and } a_{6,3} \neq 0) \text{ or } a_{6,5}= a_{6,4}a_{6,3} \neq 0
& 
& 
& \text{5A. 2.} \label{Tab5A2}
& \text{ell}
\\

& +
& 2A_1
& \text{if } a_{6,5}=a_{6,3}=0 \text{ and } a_{6,4} \neq 0
& 
& 
& \text{5A. 4.} \label{Tab5A4}
& \text{ell}
\\

& +
& A_2
& \text{if } (a_{6,5}= a_{6,4}=0 \text{ and } a_{6,3} \neq 0) \text{ or } (a_{6,3}^2 = a_{6,4} \text{ and } a_{6,3}^3= a_{6,5} \neq 0)
& 
& 
& \text{5A. 3.} \label{Tab5A3}
& \text{ell}
\\

& +
& A_3
& \text{if }  a_{6,5}= a_{6,4}= a_{6,3}=0
& 
& 
& \text{5A. 5.} \label{Tab5A5}
& \text{ell}
\\
\hline

\end{array}$$
%\vspace{-4mm}
\caption{$D_4$ and $D_5$ singularities on del Pezzo surfaces in $\Char(k)=2$} 
\label{TableD4D5-char2}
\end{table} 
\end{landscape}

\restoregeometry

\clearpage

\newgeometry{top=20 mm, bottom=20 mm, left=20 mm, right=53mm}

\begin{landscape}

\begin{table}[h] \renewcommand{\arraystretch}{1.35}
 $$
 \begin{array}{|c c c|c|c|c|c|c|} 
 \hline
\multicolumn{3}{|c|}{
\begin{array}{c}
    \text{RDP}\\
    \text{configuration}
\end{array}
}
& \begin{array}{c}
    \text{\hspace{7mm} Weierstra{\ss} equation of } X \text{ in } \bbP(1,1,2,3) \text{\hspace{7mm}} \\ 
    \hdashline
     \text{condition for extra RDPs}
\end{array}
& \Delta
& j
& \scriptstyle{\begin{array}{c}
\text{Lang's /}\\
\text{Ito's type}
\end{array} }
& \begin{array}{c}
\text{ell /}\\
\text{q-ell}
\end{array}
\\

\hline  \hline
%& 
\multicolumn{8}{|c|}{\bf{D_6}}
\\
\hline

%exception B
D_6^0
& + 
& A_1
& \multicolumn{5}{c|}{\text{occurs only in degree }$2$ \text{ (see Proposition \ref{Prop RDP iff RDP deg 1} (\hyperref[Exception B]{B.}))}}
\\
\hline

%5.2.(e)
D_6^0
& +
& 2A_1
& 
y^2 = x^3 + (t^3s+ t^2s^2)x
& 0
& 
&\text{5.2.(e)} \label{Tab52(e)}
& \text{q-ell}
\\
\hline

\hline

%13B
\multicolumn{3}{|c|}{D_6^1}
& 
y^2 + t^3y= x^3 + (a_{2,2}t^2 + a_{2,1}ts) x^2 + t^3sx
& t^{12} 
& 0
&\text{13B} \label{Tab13B-D61}
& \text{ell}
\\
\hline

\hline

%5B.
D_6^2
&
&
& 
y^2 + txy= x^3 + a_{2,1}tsx^2 + a_{6,5}t^5s + t^4s^2 \text{ with } a_{2,1} \neq 0
& 
t^{10}s (a_{6,5}t + s)
& \frac{t^{12}}{\Delta}
& \text{5B.} \label{Tab5B-D62}
& \text{ell}
\\
\hdashline

& +
& A_1
& \text{if }  a_{6,5}=0
& 
& 
& \text{5B. 2.} \label{Tab5B2}
& \text{ell}
\\

\hline  \hline
%& 
\multicolumn{8}{|c|}{\bf{D_7}}
\\
\hline

%13C
\multicolumn{3}{|c|}{D_7^1}
& 
y^2+  t^3y= x^3+ tsx^2
& t^{12} 
& 0
&\text{13C} \label{Tab13C-D71}
& \text{ell}
\\
\hline

\hline

%5C.
\multicolumn{3}{|c|}{D_7^2}
& 
y^2 + txy= x^3 + a_{2,1}tsx^2 + t^5s 
\text{ with } a_{2,1} \neq 0
& t^{11}s 
& \frac{t}{s}
&\text{5C.} \label{Tab5C-D72}
& \text{ell}
\\
\hline

\hline  \hline
%& 
\multicolumn{8}{|c|}{\bf{D_8}}
\\
\hline

%5.2.(b)
\multicolumn{3}{|c|}{D_8^0}
& 
y^2 = x^3 + t^2s^2x + t^5s
& 0
& 
&\text{5.2.(b)} \label{Tab52(b)}
& \text{q-ell}
\\
\hline

\hline

%5D.
\multicolumn{3}{|c|}{D_8^3}
& 
y^2 + t x y= x^3+ tsx^2 + a_{6,6}t^6
\text{ with } a_{6,6} \neq 0
& a_{6,6} t^{12} 
& \frac{1}{a_{6,6}}
&\text{5D.} \label{Tab5D-D83}
& \text{ell}
\\
\hline

\end{array}$$
%\vspace{-4mm}
\caption{$D_6, D_7$ and $D_8$ singularities on del Pezzo surfaces in $\Char(k)=2$} 
\label{TableD6D7D8-char2}
\end{table}

\end{landscape}

\restoregeometry

\newpage

\newgeometry{top=20 mm, bottom=20 mm, left=40 mm, right=22mm}

\begin{landscape}

\begin{table}[h] \renewcommand{\arraystretch}{1.35}
 $$
 \begin{array}{|c c c|c|c|c|c|c|} 
 \hline
\multicolumn{3}{|c|}{
\begin{array}{c}
    \text{RDP}\\
    \text{configuration}
\end{array}
}
& \begin{array}{c}
    \text{\hspace{7mm} Weierstra{\ss} equation of } X \text{ in } \bbP(1,1,2,3) \text{\hspace{7mm}} \\ 
    \hdashline
     \text{condition for extra RDPs}
\end{array}
& \Delta
& j
& \scriptstyle{\begin{array}{c}
\text{Lang's /}\\
\text{Ito's type}
\end{array} }
& \begin{array}{c}
\text{ell /}\\
\text{q-ell}
\end{array}
\\

\hline  \hline
%& 
\multicolumn{8}{|c|}{\bf{E_6}}
\\
\hline

%14
E_6^0
&
&
& 
y^2 + t^2s y = x^3 + a_{2,2}t^2 x^2 + a_{6,5}t^5s
& t^8s^4
& 0
&\text{14} \label{Tab14-E60}
& \text{ell}
\\
\hdashline

& +
& A_1
& 
\text{if } a_{6,5}= 0 \text{ and } a_{2,2} \neq 0
& 
& 
& \text{14 10A} \label{Tab1410A}
& \text{ell}
\\

& +
& A_2
& 
\text{if } a_{6,5}=  a_{2,2} = 0
& 
& 
& \text{14 11} \label{Tab1411}
& \text{ell}
\\
\hline

\hline

%6.
E_6^1
&
&
& 
y^2 + txy= x^3 + tsx^2 + a_{6,5}t^5s + a_{6,4}t^4s^2 + t^3s^3 + t^2s^4
& t^{8}s(a_{6,5}t^3 + a_{6,4}t^2s + ts^2 + s^3)
& \frac{t^{12}}{\Delta}
&\text{6.} \label{Tab6-E61}
& \text{ell}
\\
\hdashline

& +
& A_1
& 
\text{if } (a_{6,5}=0 \text{ and } a_{6,4} \neq 0) \text{ or } a_{6,5}= a_{6,4} \not\in \{0,1\}
& 
& 
& \text{6. 2.} \label{Tab62}
& \text{ell}
\\

& +
& A_2
& 
 a_{6,5}=a_{6,4} \in \{0,1\}
& 
& 
&\text{6. 3.} \label{Tab63}
& \text{ell}
\\
\hline

\hline  \hline
%& 
\multicolumn{8}{|c|}{\bf{E_7}}
\\
\hline

%exception A
\multicolumn{3}{|c|}{E_7^0}
& \multicolumn{5}{c|}{\text{occurs only in degree }$2$ \text{ (see Proposition \ref{Prop RDP iff RDP deg 1} (\hyperref[Exception A]{A.}))}}
\\
\hline

%5.2.(c)
E_7^0
& +
& A_1
& 
y^2 = x^3 + t^3s x
& 0
& 
&\text{5.2.(c)} \label{Tab52(c)}
& \text{q-ell}
\\
\hline

\hline

%15
\multicolumn{3}{|c|}{E_7^2}
& 
y^2+  t^3y= x^3+ t^3sx
& t^{12}
& 0
&\text{15} \label{Tab15-E72}
& \text{ell}
\\
\hline

\hline

%7.
E_7^3
&
&
& 
y^2 + txy= x^3 +  a_{6,5}t^5s +t^4s^2 
& t^{10}s(a_{6,5}t + s)
& \frac{t^{12}}{\Delta}
&\text{7.} \label{Tab7-E73}
& \text{ell}
\\
\hdashline

& +
&A_1
& 
\text{if } a_{6,5}=0
& 
&
&\text{7. 2.} \label{Tab72-E73}
& \text{ell}
\\
\hline

\hline  \hline
%& 
\multicolumn{8}{|c|}{\bf{E_8}}
\\
\hline

%5.2.(a)
\multicolumn{3}{|c|}{E_8^0}
& 
y^2 = x^3 + t^5s
& 0
& 
&\text{5.2.(a)} \label{Tab52(a)}
& \text{q-ell}
\\
\hline

\hline

%16
\multicolumn{3}{|c|}{E_8^3}
& 
y^2 + t^3y = x^3 + t^5s
& t^{12}
& 0
&\text{16} \label{Tab16-E83}
& \text{ell}
\\
\hline

\hline

%8.
\multicolumn{3}{|c|}{E_8^4}
& 
y^2 + t x y= x^3+ t^5s 
& t^{11}s
& \frac{t}{s}
&\text{8.} \label{Tab8-E84}
& \text{ell}
\\
\hline

\end{array}$$
%\vspace{-4mm}
\caption{$E_6, E_7$ and $E_8$ singularities on del Pezzo surfaces in $\Char(k)=2$} 
\label{TableE6E7E8-char2}
\end{table}

\end{landscape}

\newpage

\newgeometry{top=20 mm, bottom=20 mm, left=25 mm, right=35mm}
\begin{landscape}
\begin{table} \renewcommand{\arraystretch}{1.35}
 $$
 \begin{array}{|c c c|c|c|c|c|c|} 
 \hline
\multicolumn{3}{|c|}{
\begin{array}{c}
    \text{RDP}\\
    \text{configuration}
\end{array}
}
& \begin{array}{c}
    \text{\hspace{7mm} Weierstra{\ss} equation of } X \text{ in } \bbP(1,1,2,3) \text{\hspace{7mm}} \\ 
    \hdashline
     \text{condition for extra RDPs}
\end{array}
& \Delta
& j
& \scriptstyle{\begin{array}{c}
\text{Lang's /}\\
\text{Ito's type}
\end{array} }
& \begin{array}{c}
\text{ell /}\\
\text{q-ell}
\end{array}
\\

\hline  \hline
%& 
\multicolumn{8}{|c|}{\bf{E_6}}
\\
\hline

%6C
\multicolumn{3}{|c|}{E_6^0}
& y^2= x^3+t^4x+t^4s^2
& -t^{12}
& 0
& \text{6C} \label{Tab6C-E60}
& \text{ell}
\\
\hline

%6A
E_6^0
&
&
& y^2= x^3+t^3s x +a_{6,5}t^5s+t^4s^2
& -t^{9}s^3
& 0
& \text{6A} \label{Tab6A-E60}
& \text{ell}
\\
\hdashline

& +
& A_1
& \text{if }a_{6,5}=0
& 
& 
& \text{6A 5.} \label{Tab6A5}
& \text{ell}
\\
\hline

%Thm3.3. (2) in QuasiEllipticChar3
E_6^0
& +
&A_2
& y^2= x^3+ t^4s^2
& 0
& 
& \text{3.3(2)} \label{Tab33(2)}
& \text{q-ell}
\\
\hline

\hline

%6B
E_6^1
&
&
& y^2= x^3 + t^2x^2 + a_{6,5}t^5s + t^4s^2
& -t^{10}s(a_{6,5}t +s)
& \frac{t^{12}}{\Delta}
& \text{6B} \label{Tab6B-E61}
& \text{ell}
\\
\hdashline

& +
& A_1
& \text{if } a_{6,5} = 0
& 
&  
& \text{6B 2.} \label{Tab6B2}
& \text{ell}
\\
\hline

%6A
E_6^1
&
&
& 
y^2=x^3 +t^2x^2+a_{6,5}t^5s + a_{6,4}t^4s^2 + t^3s^3
& -t^9s(a_{6,5}t^2 + a_{6,4}ts +s^2)
& \frac{t^{12}}{\Delta}
& \text{6A} \label{Tab6A-E61}
& \text{ell}
\\
\hdashline

& + 
& A_1
& \text{if } a_{6,4} \neq 0 \text{ and } (a_{6,5}=0 \text{ or } a_{6,5}=a_{6,4}^2)
& 
& 
& \text{6A 2.} \label{Tab6A2}
& \text{ell}
\\

& + 
& A_2
& \text{if } a_{6,5}=0 \text{ and } a_{6,4}=0
& 
& 
& \text{6A 1.} \label{Tab6A1}
& \text{ell}
\\

\hline \hline
%& 
\multicolumn{8}{|c|}{\bf{E_7}}
\\
\hline

%7
\multicolumn{3}{|c|}{E_7^0}
& y^2 = x^3 + t^3sx + t^5s
& -t^{9}s^3
& 0
& \text{7} \label{Tab7-E70}
& \text{ell}
\\
\hline

E_7^0 
& +
& A_1
& y^2 = x^3 + t^3sx 
& -t^{9}s^3
& 0
& \text{7 5.} \label{Tab75}
& \text{ell}
\\
\hline

\hline

%7
E_7^1
&
&
& 
y^2=x^3 +t^2x^2+a_{6,5}t^5s - t^4s^2 + t^3s^3
& -t^9s(a_{6,5}t^2 -ts + s^2)
& \frac{t^{12}}{\Delta}
& \text{7} \label{Tab7-E71}
& \text{ell}
\\
\hdashline

& +
& A_1
& \text{if } a_{6,5} \in \{0,1\}
& 
& 
& \text{7 2.} \label{Tab72-E71}
& \text{ell}
\\

\hline \hline
%& 
\multicolumn{8}{|c|}{\bf{E_8}}
\\
\hline

%Thm3.3. (1) in QuasiEllipticChar3
\multicolumn{3}{|c|}{E_8^0}
& y^2= x^3+ t^5s
& 0
& 
& \text{3.3(1)}  \label{Tab33(1)}
& \text{q-ell}
\\
\hline

\hline

%8B
\multicolumn{3}{|c|}{E_8^1}
& y^2= x^3+t^4 x+ t^5s
& -t^{12}
& 0
& \text{8B} \label{Tab8B-E81}
%{\footnote{or Case I in ExtremalCharpII}}
& \text{ell}
\\
\hline

\hline

%8A
\multicolumn{3}{|c|}{E_8^2}
& y^2= x^3+t^2x^2+ t^5s
& -t^{11}s
& -\frac{t}{s}
& \text{8A 1.} \label{Tab8A1}
%{\footnote{or Case IV in ExtremalCharpII}}
& \text{ell}
\\
\hline

\end{array}$$
%\vspace{-4mm}
\caption{$E_6, E_7$ and $E_8$ singularities on del Pezzo surfaces in $\Char(k)=3$} 
\label{TableE6E7E8-char3}
\end{table} 
\end{landscape}

\newpage

\restoregeometry

\pagestyle{fancy}

\begin{table}[h] \renewcommand{\arraystretch}{1.35}
 $$
 \begin{array}{|c c c|c|c|c|c|} 
 \hline
\multicolumn{3}{|c|}{
\begin{array}{c}
    \text{RDP}\\
    \text{configuration}
\end{array}
}
& \begin{array}{c}
\text{Weierstra{\ss} equation of } X \\
\text{in } \bbP(1,1,2,3)
\end{array}
& \Delta
& j
& \scriptstyle{\begin{array}{c}
\text{Miranda's \&}\\
\text{Persson's type}
\end{array} 
% see Lang ExtremalCharpII Theorem 4.1
}
\\

\hline  \hline
%& 
\multicolumn{7}{|c|}{\bf{E_8}}
\\
\hline

\multicolumn{3}{|c|}{E_8^0}
& y^2= x^3 + t^5s
& -2t^{10}s^2
& 0
& X_{22} \label{TabX22}
\\
\hline

\hline

\multicolumn{3}{|c|}{E_8^1}
& y^2= x^3 +t^4x + t^5s
& t^{10} (t^2-2s^2)
& \frac{3t^{12}}{\Delta}
& X_{211} \label{TabX211}
\\
\hline

\end{array}$$
\vspace{-4mm}
\caption{$E_8$ singularities on del Pezzo surfaces in $\Char(k)=5$} 
\label{TableE8-char5}
\end{table}

%%%%%%%%%%%%%%%%%%%%%
% References
%%%%%%%%%%%%%%%%%%%%%

\newcommand{\etalchar}[1]{$^{#1}$}
\ifx\undefined\bysame
\newcommand{\bysame}{\leavevmode\hbox to3em{\hrulefill}\,}
\fi

\end{document}